\documentclass[mathscr]{amsart}
\usepackage{amsthm, amssymb, enumerate, eucal}

\numberwithin{equation}{section}
\theoremstyle{plain}
\newtheorem{theorem}[equation]{Theorem}
\newtheorem{lemma}[equation]{Lemma}

\newtheorem{proposition}[equation]{Proposition}
\newtheorem{corollary}[equation]{Corollary}
\theoremstyle{definition}
\newtheorem{remark}[equation]{Remark}
\newtheorem{defn}[equation]{Definition}

\newcommand{\RR}{\ensuremath{\mathbb{R}}}

\newcommand{\im}{\ensuremath{\textnormal{Im}\,}}
\newcommand{\re}{\ensuremath{\textnormal{Re}}}
\newcommand{\norm}[2]{\ensuremath{\left|\kern-.1em\left|#2\right|\kern-.1em\right|_{#1}}}

\newcommand{\paren}[1]{\ensuremath{\left( #1\right)}}
\newcommand{\abso}[1]{\ensuremath{\left| #1\right|}}
\newcommand{\braces}[1]{\ensuremath{\left\{ #1\right\}}}
\newcommand{\del}{\ensuremath{\bigtriangledown}}

\begin{document}
\title[energy-critical focusing NLS]{Global well-posedness, scattering and blow-up for the energy-critical, focusing, non-linear Schr\"odinger equation in the radial case}
\author{Carlos~E. Kenig}
\address{Department of Mathematics \\
		University of Chicago \\
		Chicago, IL 60637 \\
		USA}
\author{Frank Merle}
\address{D\'{e}partement de Math\'{e}matiques\\
		Universit\'{e} de Cergy-Pontoise\\
		Pontoise, 95302 Cergy-Pontoise\\
		FRANCE}
\thanks{The first author was supported in part by NSF, and the second one in part by CNRS. Part of this research was carried out during visits of the second author to  the University of Chicago}

\maketitle

\section{Introduction}\label{1}

In this paper, we consider the  $\dot{H}^1$ critical non-linear Schrodinger equation
\begin{equation*}\label{CP+-}
\begin{cases}
	i\partial_tu + \Delta u \pm \abso{u}^{\frac{4}{N-2}}u = 0	&(x,t)\in \RR^N\times \RR \\
	u|_{t=0} = u_0\in \dot{H}^1(\RR^N).
\end{cases}
\end{equation*}
Here the $-$ sign corresponds to the defocusing problem, while the $+$ sign corresponds to the focusing problem. The theory of the Cauchy problem (CP) for this equation was developed in \cite{CW} (Cazenave and Weissler). They show that if  $\norm{\dot{H}^1}{u_0}\le \delta$, $\delta$ small, there exists a unique solution $u \in  C(\RR;\dot{H}^1(\RR^N))$ with the norm 
$\norm{L_{x,t}^{\frac{2(N+2)}{N-2}}}{u}<\infty$(i.e. the solution scatters in $\dot{H}^1(\RR^N)$). See section 2 of this paper for a review of these results.

In the defocusing case, Bourgain (\cite{B1}, \cite{B2}) proved that, for $N=3,4$ and $u_0$ radial, this also holds for $\norm{\dot{H}^1}{u_0}<+\infty$, and that for more regular $u_0$, the solution preserves the smoothness for all time. (Another proof of this last fact is due to Grillakis \cite{Gr} for $N=3$). Bourgain's result was then extended to $N\ge5$ by Tao \cite{T}, still under the assumption that $u_0$ is radial. Then in \cite{CKSTT} (Colliander, Keel, Staffilani, Takaoka and Tao) the result was obtained for general $u_0$, when $N=3$. This was extended to $N=4$ in \cite{RV} (Ryckman, Visan) and finally to $N\ge5$ in \cite{V} (Visan).

In the focusing case, these results do not hold. In fact, the classical virial identity (see for example Glassey in \cite{G} and section 5) 
$$\frac{d^2}{dt^2}\int \abso{x}^2\abso{u_0(x,t)}^2\,dx = 8\braces{\int\abso{\del u(t)}^2 - \abso{u(t)}^{\frac{2N}{N-2}}}$$
shows that if 
$E(u_0) = \frac{1}{2} \int \abso{\del u(t)}^{2} - \frac{N-2}{2N}\int \abso{u(t)}^{\frac{2N}{N-2}}<0$ 
and $\abso{x}u_0\in L^2(\RR^N)$, the solution must break down in finite time. Moreover, 
$$W(x) = W(x,t) = \frac{1}{\paren{1+\frac{\abso{x}^2}{N(N-2)}}^{N-2/2}}$$
 is in $\dot{H}^1(\RR^N)$ and solves the elliptic equation
 $$\Delta W + \abso{W}^{\frac{4}{N-2}}W = 0,$$ so that scattering cannot always occur even for global (in time) solutions.

\medskip

In this paper we initiate the detailed study of the focusing case. We show (Corollary 5.14):

\begin{theorem}Assume that $E(u_0)<E(W)$, $\norm{\dot{H}^1}{u_0} < \norm{\dot{H}^1}{W}$,  $N=3,4,5$ and $u_0$ is radial.  Then the solution $u$  with data $u_0$ at $t=0$ is defined for all time and there exists $u_{0,+}$, $u_{0,-}$ in $\dot H^1$ such that
\[
\lim_{t\to +\infty}\norm{\dot H^1}{u(t) - e^{it\Delta}u_{0,+}} = 0,\quad \lim_{t\to -\infty}\norm{\dot H^1}{u(t) - e^{it\Delta}u_{0,-}} = 0.
\]
\end{theorem}


\medskip

Antecedents to this kind of result can be found in the $L^2$ critical case, in the work of Weinstein \cite{W} and in the $H^1$ subcritical case in the works of Beresticky and Cazenave \cite{BC}, and Zhang \cite{Z}. In particular in  \cite{BC}, the authors use variational ideas and the relationship with the virial identity.

We expect that our arguments will extend to the case of radial data, for $N\ge6$ using arguments similar to those in the appendix of \cite{T}, \cite{TV}, \cite{V} (Tao and Visan). (It remains an interesting problem to remove the radiality.) The result is optimal in that clearly the solution $W$ does not scatter. We also show that for $u_0$  radial, $E(u_0)<E(W)$, but $\norm{\dot{H}^1}{u_0} > \norm{\dot{H}^1}{W},$ the solution must break down in finite time.

Our proof introduces a new point of view for these problems. Using a concentration compactness argument (section 4), we reduce matters to a rigidity theorem, which we prove in section 5, with the aid of a localized virial identity (in the spirit of Merle \cite{M}, \cite{M1}). The radiality enters only at one point, in our proof of the rigidity theorem (see Remark 5.2). 
We think that the general strategy of our proof with one extra ingredient should also apply in the non-radial case. In section 3, we prove some elementary variational estimates which yield the necessary coercivity for our arguments. These are automatic in the defocusing case and thus our proof gives an alternative approach to \cite{B1} and \cite{T} for $N=3,4,5$.

Acknowledgment: We would like to thank the referees for their suggestions and their careful reading of the manuscript.

\section{A review of the Cauchy problem}\label{2}
In this section we will review the Cauchy problem
\begin{equation*}\label{CP}
\tag{CP}
\begin{cases}
	i\partial_tu + \Delta u + \abso{u}^{\frac{4}{N-2}}u = 0	&(x,t)\in \RR^N\times \RR \\
	u|_{t=0} = u_0\in \dot{H}^1(\RR^N)
\end{cases}
\end{equation*}
i.e., the $\dot{H}^1$ critical, focusing, Cauchy problem for NLS.  We need two preliminary results.
\begin{lemma}[Strichartz estimate \cite{C},\cite{KT}]\label{2.1}
We say that a pair of exponents $(q,r)$ is admissible if $\frac{2}{q}+\frac{N}{r} = \frac{N}{2}$ and $2\le q$, $r\le \infty$.  Then, if $2\le r\le \frac{2N}{N-2}$ ($N\ge 3$) (or $2\le r<\infty$, $N=2$ and $2\le r\le \infty$, $N=1$) we have
\begin{enumerate}[i)]
\item 
\[
	\norm{L_t^qL_x^r}{e^{it\Delta}h}\le C\norm{L^2}{h}
\]
\item
\[
	\norm{L_t^qL_x^r}{\int_{-\infty}^{+\infty} e^{i(t-\tau)\Delta}g(-,\tau)\,d\tau} +
	\norm{L_t^qL_x^r}{\int_{0}^{t} e^{i(t-\tau)\Delta}g(-,\tau)\,d\tau}\le C\norm{L_t^{q'}L_x^{r'}}{g}
\]
\item
\[
	\norm{L_x^2}{\int_{-\infty}^{+\infty} e^{it\Delta}g(-,\tau)\,d\tau}\le C\norm{L_t^{q'}L_x^{r'}}{g}.
\]
\end{enumerate}
\end{lemma}

\begin{lemma}[Sobolev Embedding]\label{2.2}
For $v\in C_0^\infty(\RR^{N+1})$, we have
\[
	\norm{L_t^{\frac{2(N+2)}{N-2}}L_x^{\frac{2(N+2)}{N-2}}}{v}\le C\norm{L_t^\frac{2(N+2)}{N-2}L_x^{\frac{2N(N+2)}{N^2+4}}}{\del_xv}\quad(N\ge 3).
\]
Note that $\frac{2(N+2)}{N-2}=q$, $\frac{2N(N+2)}{N^2+4}=r$ is admissible.)
\end{lemma}

\begin{remark}\label{2.3}
Let $f(u) = \abso{u}^{\frac{4}{N-2}}u$, then clearly \\
$\abso{f(u)}\le \abso{u}^{N+2/N-2}$, $\abso{\partial_zf(u)}\le C\abso{u}^{\frac{4}{N-2}}$, $\abso{\partial_{\overline{z}}f(u)}\le C\abso{u}^{\frac{4}{N-2}}$. \\
Moreover, for $3\le N\le 6$,
\[
	\left. \begin{matrix}\abso{\partial_zf(u) - \partial_zf(v)}\\ \abso{\partial_{\overline{z}}f(u)-\partial_{\overline{z}}f(v)}\end{matrix}\right\}\le C \abso{u-v}\cdot \left\{\abso{u}^{\frac{6-N}{N-2}}+\abso{v}^{\frac{6-N}{N-2}}\right\}.
\]
Also, note that $(\del f)(u(x)) = \partial_zf(u(x))u(x) + \partial_{\overline{z}}f(u(x))\overline{u}(x)$, so that $\abso{f(u)-f(v)}\le \abso{u-v}\braces{\abso{u}^{\frac{4}{N-2}}+\abso{v}^{\frac{4}{N-2}}}$.  Moreover, \\
\[\begin{aligned}
\del_x(f(u(x)))&-\del_x(f(v(x))) = 
(\del f)(u(x))\del u - (\del f)(v(x))\del v \\
&= (\del f)(u(x))\del u - (\del f)(u(x))\del v + \braces{\del f(u(x))) - \del f(v(x))}\del v,
\end{aligned}\]
so 
$\abso{\del_x f(u(x)) - \del_xf(v(x))}\le C\abso{u(x)}^{\frac{4}{N-2}}\abso{\del u-\del v} \\ 
+C\abso{\del v}\braces{\abso{u}^{\frac{6-N}{N-2}} + \abso{v}^{\frac{6-N}{N-2}}}\abso{u-v}$.
\end{remark}

\begin{remark}\label{2.4}
In the estimate ii) in Lemma 2.1, one can actually show: (\cite{KT})
\begin{enumerate}[ii')]
\item
\[
	\norm{L_t^qL_x^r}{\int_{-\infty}^{+\infty}e^{i(t-\tau)\Delta} g(-,\tau)\,d\tau}\le C\norm{L_t^{m'}L_x^{n'}}{g},
\]
\end{enumerate}
where $(q,r)$, $(m,n)$ are any pair of admissible indices as in i) of Lemma 2.1.
\end{remark}

Let us define $S(I), W(I)$ norm for an interval $I$ by
$$ \norm{S(I)}{v} = \norm{L_I^{\frac{2(N+2)}{N-2}}L_x^{\frac{2(N+2)}{N-2}}}{v}     \  \  \mbox{and}  \   \
\norm{W(I)}{v} = \norm{L_I^{\frac{2(N+2)}{N-2}}L_x^{\frac{2N(N+2)}{N^2+4}}}{v}.$$

\begin{theorem}[See \cite{CW}]
Assume $u_0\in \dot{H}^1(\RR^N)$, $t_0\in I$ an interval, and $\norm{\dot{H}^1}{u_0}\le A$.  Then, (for $3\le N\le 5$) there exists $\delta = \delta(A)$ such that, if 
$\norm{S(I)}{e^{i(t-t_0)\Delta} u_0}
<\delta$, there exists a unique solution $u$ to (CP) in $I\times \RR^N$, with $u\in C(I;\dot{H}^1(\RR^N))$, 
$$
\norm{W(I)}{\del_xu} 
<\infty, \ \ \  
\norm{S(I)}{u}
\le 2\delta.$$
  Moreover, if $u_{0,k}\to u_0$ in $\dot{H}^1$ (so that, as we will see, for $k$ large $\norm{S(I)}{e^{i(t-t_0)\Delta}u_{0,k}}<\delta$) the corresponding solutions $u_k\to u$ in $C(I;\dot{H}^1(\RR^N))$.
\end{theorem}
\begin{proof}[Sketch of Proof]
Let us assume, without loss of generality, that $t_0=0$.  (CP) is equivalent to the integral equation
\[
	u(t) = e^{it\Delta}u_0 + \int_0^t e^{i(t-t')\Delta}f(u)\,dt',
\]
where $f(u) = \abso{u}^{\frac{4}{N-2}}u$.  We now let 
$B_{a,b} = \left\{v\text{ on }I\times \RR^n~:~ \norm{S(I)}{v}\le a,\right.$\\ $\left.\norm{W(I)}{\del v}\le b\right\}$
and $\Phi_{u_0}(v) = e^{it\Delta}u_0 + \int_0^t e^{i(t-t')\Delta}f(v)\,dt'$.  
We will next choose $\delta,a,b$ so that $\Phi_{u_0}(v):B_{a,b}\to B_{a,b}$ and is a contraction there:  note that $$\norm{W(I)}{\del \Phi_{u_0}(v)}\le CA + C\norm{L_I^2L_x^{\frac{2N}{N+2}}}{\del_xf(v)}.$$  This follows, for the first term, by i) ($q = \frac{2(N+2)}{N-2}$, $r = \frac{2N(N+2)}{N^2+4}$) in Lemma 2.1 and by ii) in Remark 2.4, with the same $q,r$ and $m' = 2$, $n' = \frac{2N}{N+2}$.  But $\del_x f(u(x)) = (\del f)(u(x))\del_xu = O(\abso{\del u}\cdot \abso{u}^{\frac{4}{N-2}})$ so that, using H\"older inequality we obtain: (for $v\in B_{a,b}$)
$$\norm{W(I)}{\del \Phi_{u_0}(v)}
	\le CA + C\norm{S(I)}{v}^{\frac{4}{N-2}}\cdot \norm{W(I)}{\del v} \le CA + Ca^{\frac{4}{N-2}}b.$$
Using Lemma 2.2 for the second term in $\Phi_{u_0}$, and the argument above together with our assumption on $u_0$ for the first term, we obtain:
\[
	\norm{S(I)}{\Phi_{u_0}(v)}\le \delta + Ca^{\frac{4}{N-2}}b.
\]
Now choose $b = 2AC$, and $a$ so that $Ca^{\frac{4}{N-2}}\le 1/2$.  Then $\norm{W(I)}{\del \Phi_{u_0}(v)}\le b$.  Next, if $\delta = a/2$, and $Ca^{(\frac{4}{N-2}-1)} b\le 1/2$ {\it (possible if $N<6$) } we obtain $\norm{S(I)}{\Phi_{u_0}(v)}\le a$, so that $\Phi_{u_0}:B_{a,b}\to B_{a,b}$.  Next, for the contraction, we use the same argument in conjunction with Remark 2.3.
\[\begin{aligned}
	&\norm{W(I)}{\del \Phi_{u_0}(v) - \del \Phi_{u_0}(v')}
	\le C\norm{L_I^2L_x^{\frac{2N}{N+2}}}{\del_xf(v)-\del_xf(v')}\\
	&\le C\norm{L_I^2L_x^{\frac{2N}{N+2}}}{\abso{v}^{\frac{4}{N-2}}\abso{\del v-\del v'}} 
	+C\norm{L_I^2L_x^{\frac{2N}{N+2}}}{\abso{v-v'}\abso{v}^{\frac{6-N}{N-2}}\abso{\del v'}} \\
	&+ C\norm{L_I^2L_x^{\frac{2N}{N+2}}}{\abso{v-v'}\abso{v'}^{\frac{6-N}{N-2}}\abso{\del v'}}.
\end{aligned}\]
The first term is bounded as before by $C\norm{S(I)}{v}^{\frac{4}{N-2}}\norm{W(I)}{\del v-\del v'}$.  For the second and third terms we use H\"older's inequality to bound them by \\ $C\norm{S(I)}{v-v'}\norm{W(I)}{\del v'} (\norm{S(I)}{v}^{\frac{6-N}{N-2}}+\norm{S(I)}{v'}^{\frac{6-N}{N-2}})$ so that
\[\begin{aligned}
	\norm{W(I)}{\del \Phi_{u_0}(v) - \del \Phi_{u_0}(v')}
	&\le Ca^{\frac{4}{N-2}}\norm{W(I)}{\del v-\del v'}  \\
	&+ Ca^{\frac{6-N}{N-2}} b\norm{S(I)}{v-v'}.
\end{aligned}\]
Lemma 2.2 gives
\[\begin{aligned}
	\norm{S(I)}{\Phi_{u_0}(v) - \Phi_{u_0}(v')}
	&\le C\norm{W(I)}{\del\Phi_{u_0}(v) - \del\Phi_{u_0}(v')}  \\
	&\le Ca^{\frac{4}{N-2}}\norm{W(I)}{\del v-\del v'} + Ca^{\frac{6-N}{N-2}} b\norm{S(I)}{v-v'}
\end{aligned}\]
and thus we establish the contraction property ($N<6$).  We then find $u\in B_{a,b}$ solving $\Phi_{u_0}(u)=u$.  To show that $u\in C(I;\dot{H}^1)$, note that $e^{it\Delta}u_0\in C(I;\dot{H}^1)$ with norm bounded by $A$.  For the term $\int_0^t e^{i(t-t')\Delta}f(u)\,dt'$, we use iii) in Lemma 2.1, with $(q',r') = (2,2N/N+2)$.  The proof of Theorem 2.5 is easily concluded from this.
(The last continuity statement is an easy consequence of the fixed point argument, see also Remark 2.17.)
\end{proof}

\begin{remark}\label{2.6}
Using Remark 2.4, it is easy to see that $\del u\in L_I^qL_x^r$ for any admissible index pair $(q,r)$.
\end{remark}

\begin{remark}\label{2.7}
There exists $\widetilde{\delta}$ such that if $\norm{\dot{H}^1}{u_0}\le \widetilde{\delta}$, the conclusion of Theorem 2.5 applies to any interval $I$.  In fact, $\norm{S(I)}{e^{it\Delta}u_0}\le C\norm{W(I)}{\del e^{it\Delta}u_0}\le C\widetilde{\delta}$, by virtue of Lemma 2.1 i) and the claim follows.
\end{remark}

\begin{remark}\label{2.8}
Given $u_0\in \dot{H}^1$, there exists $(0\in)~I$ such that the hypotheses  of Theorem 2.5 is verified on $I$.  This is clear because of $\norm{S(I)}{e^{it\Delta}u_0}\le C\norm{W(I)}{\del e^{it\Delta}u_0}$ and the fact that $\norm{W(\RR)}{\del e^{it\Delta}u_0} < \infty$ by Lemma 2.1 i).
\end{remark}

\begin{remark}[Energy Identity]\label{2.9}
If $u$ is the solution constructed in Theorem 2.5, we have (with $\frac{1}{2^*} = \frac{1}{2}-\frac{1}{N}$) that 
$$E(u(t)) = \int_{\RR^N}\braces{\frac{1}{2}\abso{\del u(t,x)}^2 - \frac{1}{2^*}\abso{u(t,x)}^{2^*}}\,dx$$
 is constant for $t\in I$.  If $u_0\in C_0^\infty(\RR^N)$ this follows from a classical integration by parts, the general general case follows from a limiting argument.
\end{remark}

\begin{defn}\label{2.10}
Let $t_0\in I$.  We say that $u\in C(I;\dot{H}^1(\RR^N))\cap \braces{\del u\in W(I)}$ is a solution of the (CP) if 
$$u|_{t_0} = u_0, \mbox{ \ \ and \ \ } u(t) = e^{i(t-t_0)\Delta}u_0 + \int_{t_0}^te^{i(t-t')\Delta}f(u)\,dt',$$
 with $f(u) = \abso{u}^{\frac{4}{N-2}}u$.  Note that if $u^{(1)}$, $u^{(2)}$ are solutions of (CP) on $I$, $u^{(1)}(t_0) = u^{(2)}(t_0)$, then $u^{(1)}\equiv u^{(2)}$ on $I\times \RR^N$.  This is because we can partition $I$ into a finite collection of subintervals $I_j$, so that, with $A = \sup_{t\in I}\max_{i=1,2}\norm{\dot{H}^1}{u^{(i)}(t)}$, the $S(I_j)$ norm of $u^{(i)}$ and the $W(I_j)$ norm of $\del u^{(i)}$ are less than $a,b$, where $a,b$ are obtained in the proof of Theorem 2.5.  If $j_0$ is then such that $t_0\in I_{j_0}$, the uniqueness of the fixed point in the proof of Theorem 2.5, combined with Remark 2.8 gives an interval $\widetilde{I}\ni t_0$ so that $u^{(1)}(t) = u^{(2)}(t)$, $t\in \widetilde{I}$.  A continuation argument now easily gives $u^{(1)}\equiv u^{(2)}$, $t\in I$.  This allows us to define a maximal interval $I(u_0) = (t_0-T_-(u_0),t_0+T_+(u_0))$, with $T_\pm(u_0)>0$, where the solution is defined.  If $T_1<t_0+T_+(u_0)$, $T_2 > t_0 - T_-(u_0)$, $T_2<t_0<T_1$, then $u$ solves (CP) in $[T_2,T_1]\times \RR^N$, so that $u\in C([T_2,T_1]\times \RR^N)$, $\del u\in W([T_2,T_1])$ and $u\in S([T_2,T_1])$.
\end{defn}

\begin{lemma}[Standard finite blow-up criterion, see \cite{C}]\label{2.11}
If $T_+(u_0)<\infty$, then  
$$\norm{S([t_0,t_0+T_+(u_0)])}{u} = +\infty. $$ 
A corresponding result holds for $T_-(u_0)$.
\end{lemma}
\begin{proof}[Sketch of Proof]
Assume $T_+(u_0)<+\infty$ and that $\norm{S([t_0,t_0+T_+(u_0)])}{u} < +\infty$.  Let $M = \norm{S([t_0,t_0+T_+(u_0)])}{u}$ and, for $\epsilon$ to be chosen, find $N = N(\epsilon)$ intervals $I_j$, $\bigcup_{j=1}^N I_j = [t_0,t_0+T_+(u_0)]$, such that $\norm{S(I_j)}{u}\le \epsilon$.  Our first step is to show that $\norm{L^\infty([t_0,t_0+T_+(u_0)];\dot{H}^1)}{u} + \norm{W([t_0,t_0+T_+(u_0)])}{\del u} < \infty$.  We write the integral equation on each interval $I_j$, to deduce (using the proof of Theorem 2.5 and iii) in Lemma 2.1) that
\[
	\sup_{t\in I_j} \norm{\dot{H}^1}{u(t)} + \norm{W(I_j)}{\del u}\le C\norm{\dot{H}^1}{u(t_j)} + C\norm{S(I_j)}{u}^{\frac{4}{N-2}}\cdot \norm{W(I_j)}{\del u},
\]
where $t_j$ is any fixed point in $I_j$.  Our desired estimate follows inductively then, by choosing $C\epsilon^{\frac{4}{N-2}}\le 1/2$.  Once the first step is done, we then choose $t_n\uparrow t_0+T_+(u_0)$ and show, using the integral equation once more, that \\ $\norm{S([t_n,t_0+T_+(u_0)])}{e^{i(t-t_n)\Delta}u(t_n)}\le \delta/2$, for $n$ large.  But then, for $n$ large but fixed, and some $\epsilon_0 > 0$, $\norm{S([t_n,t_0+T_+(u_0)]+\epsilon_0)}{e^{i(t-t_n)\Delta}u(t_n)}\le \delta$.  Now, Theorem 2.5 applies and together with Definition 2.10 we reach a contradiction.
\end{proof}

\begin{defn}\label{2.12}
Let $v_0\in \dot{H}^1$, $v(x,t) = e^{it\Delta}v_0$ and let $\{t_n\}$ be a sequence, with $\lim_{n\to \infty}t_n = \overline{t}\in [-\infty,+\infty]$.  We say that $u(x,t)$ is a non-linear profile associated with $(v_0,\{t_n\})$ if there exists an interval $I$, with $\overline{t}\in I$ (if $\overline{t} = \pm \infty$, $I = [a,+\infty)$ or $(-\infty, a]$) such that $u$ is a solution of (CP) in $I$ and
\[
	\lim_{n\to \infty}\norm{\dot{H}^1}{u(-,t_n)-v(-,t_n)} = 0.
\]
\end{defn}

\begin{remark}\label{2.13}
There always exists a non-linear profile associated to $(v_0,\{t_n\})$.  In fact, if $\overline{t}\in (-\infty,+\infty)$, this is clear by Remark 2.8, with $u_0 = v(x,\overline{t})$.  If $\overline{t} = +\infty$, we solve the integral equation
\[
	u(t) = e^{it\Delta}v_0 + \int_t^{+\infty}e^{i(t-t')\Delta}f(u)\,dt'
\]
in $(t_{n_0},+\infty)\times \RR^N$, for $n_0$ so large that $\norm{S((t_{n_0},\infty))}{e^{it\Delta}v_0}\le \delta$, where $\delta$ is as in Theorem 2.5.  Then, if $n$ is large $u(t_n) - v(t_n) = \int_{t_n}^{+\infty}e^{i(t-t')\Delta}f(u)\,dt'$, and we have $\norm{L^2_{(t>t_{n_0})}L_x^{2N/N+2}}{f(u)}<\infty$, as in the proof of Theorem 2.5.  But then, using iii) in Lemma 2.1 we obtain $\norm{\dot{H}^1}{u(t_n) - v(t_n)}\le C\norm{L^2_{(t>t_n)}L_x^{2N/N+2}}{f(u)}$, which clearly goes to $0$ as $n$ goes infinity.  A similar argument applies when $\overline{t} = -\infty$.

Note also that if $u^{(1)}$, $u^{(2)}$ are both non-linear profiles associated to $(v_0,\{t_n\})$ in an interval $I$ with $\overline{t}\in I$, then $u^{(1)}\equiv u^{(2)}$ on $I$.  In fact, if $\overline{t}\in (-\infty, +\infty)$, this is clear from the Definition 2.13 and the uniqueness result in Definition 2.10.  If $\overline{t} = +\infty$, since $\norm{W(I)}{\del u^{(i)}}<\infty$, for $n\ge n_0$, we have $\norm{W(t_n,+\infty)}{\del u^{(i)}}\le \widetilde{\delta}$, where $\widetilde{\delta}$ is as small as we like.  By the proof of Theorem 2.5, we have (with a constant independent of $u$) that for $n\gg n_0$
\[
	\sup_{t\in (t_{n_0},t_n)}\norm{L^2}{\del u^{(1)}(t) - \del u^{(2)}(t)}\le C\norm{L^2}{\del u^{(1)}(t_n) - \del u^{(2)}(t_n)}.
\]
This easily shows that $u^{(1)}\equiv u^{(2)}$ on $(t_{n_0},+\infty)$ and hence on $I$, as claimed.  The case $\overline{t} = -\infty$ is similar.  Because of this remark, we can always define a maximal interval $I$ of existence for the non-linear profile associated to $(v_0,\{t_n\})$.  If $\overline{t}\in (-\infty,+\infty)$, $I = (a,b)$, $I'\Subset I$, then $\sup_{t\in I'}\norm{\dot H^1}{u(t)}<\infty$, $\norm{S(I')}{u}<\infty$, $\norm{W(I')}{\del u}<\infty$, but if either $a$ or $b$ are finite $\norm{S(I)}{u} = +\infty$.  If $\overline{t} = \pm\infty$, say $\overline{t} = +\infty$, $I = (a,+\infty)$, $I' = (\alpha,+\infty)$, $\alpha>a$, similar statements can be made.  If $a>-\infty$, we can also say $\norm{S(I)}{u} = +\infty$.
\end{remark}

\begin{theorem}[Long-time perturbation theory, see also \cite{TV}]\label{2.14}
Let $I\subset \RR$ be a time interval and let $t_0\in I$.  Let $\widetilde{u}$ be defined on $I\times \RR^N$ ($3\le N\le 5$) and satisfy $\sup_{t\in I}\norm{\dot H^1}{\widetilde{u}}\le A$, $\norm{S(I)}{\widetilde{u}}\le M$ for some constants $M,A>0$.  Assume that
\[
	(i\partial_t\widetilde{u} + \Delta \widetilde{u} + f(\widetilde u)) = e\qquad (t,x)\in I\times \RR^N
\]
(in the sense of the appropriate integral equation) and that 
$$\norm{\dot H^1}{u_0-\widetilde u(t_0)}\le A', \ \  \norm{L_I^2L_x^{\frac{2N}{N+2}}}{\del e}\le \epsilon, \ \  \norm{S(I)}{e^{i(t-t_0)\Delta}[u_0-\widetilde u(t_0)]}\le \epsilon.$$  
Then, there exists \\ $\epsilon_0 = \epsilon_0(M,A,A',N)$ such that there exists a solution of (CP) with $u(t_0) = u_0$ in $I\times \RR^N$, for $0<\epsilon<\epsilon_0$, with $\norm{S(I)}{u}\le C(M,A,A',N)$ and 
$\forall t \in I, \norm{\dot H^1}{u(t)-\widetilde u(t)}\le C(A,A',M,N) (A'+ \epsilon)$.
\end{theorem}
\begin{proof}
We start the proof by showing that $\norm{W(I)}{\del \widetilde{u}}\le \widetilde M$, where $\widetilde M = \widetilde M(A,M,N)$, for $\epsilon\le \epsilon_0$.  In fact, for $\eta = \eta(N)$ small, to be determined, split $I$ into $\gamma = \gamma(M,\eta)$ interval $I_j$ so that $\norm{S(I_j)}{\widetilde u}\le \eta$.  Using the integral equation, we have
\[
	\norm{W(I_j)}{\del \widetilde u}\le A+C\norm{S(I_j)}{\widetilde u}^{\frac{4}{N-2}}\norm{W(I_j)}{\del \widetilde u}+C\norm{L_I^2L_x^{2N/N+2}}{\del e},
\]
as in the proof of Theorem 2.5, and the claim follows if $C\eta^{\frac{4}{N-2}} < 1/2$.  Next, we write $u = \widetilde u + w$ and notice that
\[
	i\partial_t w + \Delta w + [f(\widetilde u+w)-f(\widetilde u)] = e.
\]
Let $I_j = [a_j,a_{j+1}]$, so that, in order to solve for $w$ we need to solve, in $I_j$, the integral equation
\[
	w(t) = e^{i(t-a_j)\Delta}w(a_j) + \int_{a_j}^te^{i(t-t')\Delta}[ f(\widetilde u + w) - f(\widetilde u)]\,dt' + \int_{a_j}^t e^{i(t-t')\Delta}e\,dt'.
\]
The proof of Theorem 2.5 {\it (which holds for $3 \le N \le 5$)} now shows that, for $\eta = \eta(N)$ small enough, and $\epsilon_0 = \epsilon_0(N)$ small enough, we can solve the integral equation (assuming $t_0=a_1$ say) in $I_1$ and obtain $w$ with the bounds $\norm{S(I)}{w}\le 2\epsilon$, $\norm{W(I_1)}{\del w}\le C(A,A')$, $\sup_{t\in I_1}\norm{\dot H^1}{w(t)}\le C(A,A')(A'+ \epsilon)$.  We now estimate $\norm{S(I_2)}{e^{i(t-a_2)\Delta}w(a_2)}$, using the integral equation.  Since $e^{i(t-a_2)\Delta}e^{i(a_2-t_0)\Delta}w(t_0) = e^{i(t-t_0)\Delta}w(t_0)$, by assumption $\norm{S(I_2)}{e^{i(t-a_2)\Delta}e^{i(a_2-t_0)\Delta}w(t_0)}\le \epsilon$.  For the integral term, we use Lemma 2.1, iii) to obtain a bound for its $\dot H^1$ norm at $a_2$ by $C\norm{S(I_1)}{w}^{\frac{4}{N-2}} \norm{S(I_1)}{\del w}\le C(2\epsilon)^{\frac{4}{N-2}} C(A,A')$.  Clearly this procedure can be iterated $\gamma = \gamma(M,N)$ times, provided $\epsilon_0$ is small enough, yielding the theorem.
\end{proof}

\begin{remark}[See \cite{C}]\label{2.15}
If $u$ is a solution of (CP) in $I\times \RR^N$, $I = [a,+\infty)$ (or $I = (-\infty,a]$) there exists $u^+\in \dot H^1$ such that 
$$\lim_{t\to +\infty}\norm{\dot H^1}{u(t) - e^{it\Delta}u_+} = 0.$$
To see this, note that $\del f(u)\in W(I)$ and hence $\norm{\dot H^1}{\int_t^\infty e^{i(t-t')\Delta}f(u)\,dt'}\to 0$ as $t\to +\infty$.  Then, $u(t) = e^{i(t-a)\Delta}u_0 + \int_a^t e^{i(t-t')\Delta}f(u)\,dt'$ and hence $u^+ = e^{-ia\Delta}u_0 + \int_a^\infty e^{-it'\Delta}f(u)\,dt'$ has the desired property.  In fact note that the argument used at the beginning of the proof of Theorem 2.14 shows that it suffices to assume $u$ to be a solution of (CP) in $I'\times \RR^N$, $I'\Subset I$, such that $\norm{S(I)}{u}<\infty$.
\end{remark}

\begin{remark}\label{2.16}
We recall that, since we are working in the focusing case, we have from the argument of Glassey \cite{G} that if $\int |x|^2|u_0|^2 < + \infty$, $E(u_0)< 0$, there exists a finite time $T$ such that the solution cannot be extended for $t>T$. Clearly, for such a $u_0$, the maximal interval of existence must be finite.  (See Definition 2.10.) Note that it is unknown if $\lim_{t\uparrow T}\norm{\dot H^1}{u(t)} = +\infty$ for a general initial data that doesn't exist for all time.
\end{remark}

\begin{remark}\label{2.17}
Theorem 2.14 also yields the following continuity fact, which will be used later: let $\widetilde{u}_0  \in \dot H^1$, $\norm{\dot H^1}{\widetilde{u}_0}\le A$, and let $\widetilde{u}$ be the solution of (CP), with maximal inteval of existence $(T_-(\widetilde{u}_0),T_+(\widetilde{u}_0))$ (see definition 2.10). Let 
$u_{0,n}\rightarrow\widetilde{u}_0$ in $\dot H^1$, 
and let $u_n$ be the corresponding solution of (CP), with maximal interval of existence $(T_-(u_{0,n}),T_+(u_{0,n}))$. Then, $T_-(\widetilde{u}_0) \le \varliminf_{n\to + \infty}T_-(u_{0,n})$, $T_+(\widetilde{u}_0) \le \varliminf_{n\to + \infty}T_+(u_{0,n})$ and for each $t \in (T_-(\widetilde{u}_0),T_+(\widetilde{u}_0)),$ $u_{n}(t)\rightarrow\widetilde{u}(t)$ in $\dot H^1$. \\
Indeed, Let 
$I \subset\subset(T_-(\widetilde{u}_0),T_+(\widetilde{u}_0))$, 
so that 
$Sup_{t \in I}\norm{\dot H^1}{\widetilde u(t)}= A < + \infty, \norm{S(I)}{\widetilde u(t)}= M < + \infty.$
We will show that, for n large, $u_n$ exists on $I$ and 
$\forall t \in I, \norm{\dot H^1}{u_n(t)-\widetilde u(t)}\le C(M,A,N) \norm{\dot H^1}{u_{0,n}-\widetilde u_0}$. 
This clearly yields the remark. To show this, apply Theorem 2.14, with $u=u_n,$ $u_0=u_{0,n}$. Then, if 
$\epsilon_0 = \epsilon_0(M,A,2A,N)$ 
and $n$ is so large that 
$\norm{\dot H^1}{u_{0,n}-\widetilde u_0} \le \epsilon_0$ 
and 
$\norm{S(I)}{e^{it\Delta}[u_{0,n}-\widetilde u_0]}\le \epsilon_0$, using the uniqueness of the solutions we obtained in Definition 2.10, the claim follows.

\end{remark}

\section{Some variational estimates}\label{3}
Let $$W(x) = W(x,t) = \frac{1}{\paren{1+\frac{\abso{x}^2}{N(N-2)}}^{(N-2)/2}},$$ 
be a stationary solution of (CP).  That is, $W$ solves the non-linear elliptic equation
\begin{equation}\label{3.3}
\Delta W + \abso{W}^{\frac{4}{N-2}}W = 0.
\end{equation}
Moreover, $W\ge 0$ and it is radially symmetric and decreasing. Note that $W\in \dot H^1$, but $W$ need not belong to $L^2(\RR^N)$.  By invariances of the equation, for $\theta_0\in [-\pi,\pi]$, $\lambda_0>0$, $x_0\in \RR^N$, $W_{\theta_0,x_0,\lambda_0}(x) = e^{i\theta_0}\lambda_0^{(N-2)/2}W(\lambda_0(x-x_0))$ is still a solution.  By the work of Aubin \cite{A}, Talenti \cite{Ta} we have the following characterization of $W$:
\begin{equation}\label{3.1}
\forall u\in \dot H^1,\quad \norm{L^{2^*}}{u}\le C_N\norm{L^2}{\del u};
\end{equation}
moreover,
\begin{equation}\label{3.2}
\text{If } \norm{L^{2^*}}{u} = C_N \norm{L^2}{\del u},~u\ne 0,\text{ then }\exists (\theta_0,\lambda_0,x_0)\text{ such that } u = W_{\theta_0,x_0,\lambda_0}, 
\end{equation}
where $C_N$ is the best constant of the Sobolev inequality in dimension $N$.

The equation (3.1) gives $\int \abso{\del W}^2 = \int\abso{W}^{2^*}$.  Also, (3.3) yields $C_N^2\int\abso{\del W}^2 = \paren{\int\abso{W}^{2^*}}^{N-2/N}$, so that $C_N^2\int\abso{\del W}^2 = \paren{\int\abso{\del W}^2}^{\frac{N-2}N}$.  Hence, 
$$\int \abso{\del W}^2 = \frac1{C_N^N} \mbox{ \ \  and \ \ } E(W) = \paren{\frac{1}{2} - \frac{1}{2^*}}\int\abso{\del W}^2 = \frac{1}{N}\frac{1}{C_N^N}.$$

\begin{lemma}\label{3.4}
Assume that $$ \norm{L^2}{\del u}^2 <  \norm{L^2}{\del W}^2.$$  
Assume moreover that $E(u)\le (1-\delta_0)E(W)$ where $\delta_0>0$.  Then, there exists $\overline{\delta} = \overline{\delta}(\delta_0,N)>0$ such that
\begin{equation}\label{3.5}
\int\abso{\del u}^2\le (1-\overline{\delta})\int\abso{\del W}^2
\end{equation}
and
\begin{equation}\label{3.6}
\int\abso{\del u}^2 - \abso{u}^{2^*}\ge \overline{\delta}\int\abso{\del u}^2
\end{equation}
and
\begin{equation}\label{3.7}
E(u)\ge 0.
\end{equation}
\end{lemma}
\begin{proof}

Consider the function $f_1(y) = \frac{1}{2}y - \frac{C_N^{2^*}}{2^*}y^{\frac{2^*}2}$, and let $\overline{y} = \norm{L^2}{\del u}^2$.  Because of (3.2), $f_1(\overline{y})\le E(u)\le (1-\delta_0)E(W) = (1-\delta_0)\frac{1}{N} \frac{1}{C_N^N}$.  Note that $f_1(0) = 0$, $f_1'(y) = \frac{1}{2} - \frac{C_N^{2^*}}{2}y^{\frac{2^*}{2}-1}$, so that $f_1'(y) = 0$ if and only if $y = y_C$, where $y_C = \frac{1}{C_N^N} = \int\abso{\del W}^2$.  Note also that $f_1(y_C) = \frac{1}{NC_N^N} = E(W)$.  But then, since $0<\overline{y}<y_C$ and $f_1(\overline{y})\le (1-\delta_0)f_1(y_C)$ and $f_1$ is nonnegative and strictly increasing between 0 and $y_C$, $f_1''(y_C)\ne 0$, we have $0<f_1(\overline{y})$ and $\overline{y}\le (1-\overline{\delta})\int\abso{\del W}^2$.  Thus (3.5) and (3.7) hold.  To show (3.6), consider the function $g_1(y) = y-C_N^{2^*}y^{N/N-2}$.  Because of (3.2) we have that $\int\abso{\del u}^2-\abso{u}^{2^*}\ge \int \abso{\del u}^2 - C_N^{2^*}\paren{\int\abso{\del u}^2}^{2^*/2} = g_1(\overline{y})$.  Note that $g_1(y) = 0$ if and only if $y=0$ or $y = y_C$ and that $g_1'(0) = 1$, $g_1'(y_C) = -\frac2{N-2}$.  We then have, for $0<y<y_C$, $g_1(y)\ge C\min\{y,(y_C-y)\}$, and so, since $0 \le \overline{y}<(1-\overline{\delta})y_C$ by (3.5), (3.6) follows. Note that $\overline{\delta} \simeq {\delta_0}^{\frac 12}$.
\end{proof}

Note that the relevance of (3.6) comes from the virial identity (see introduction).

\begin{corollary}\label{3.8}
Assume that $u\in \dot H^1$ and that $\int\abso{\del u}^2 < \int\abso{\del W}^2$.  Then $E(u)\ge 0$.
\end{corollary}
\begin{proof}
If $E(u)\ge E(W) = \frac{1}{NC_N^N}$, this is obvious.  If $E(u) < E(W)$, the claim follows from (3.7).
\end{proof}

\begin{theorem}[Energy trapping]\label{3.9}
Let $u$ be a solution of the (CP), with $t_0 = 0$, $u|_{t=0} = u_0$ such that for  $\delta_0>0$
$$\int\abso{\del u_0}^2 < \int\abso{\del W}^2 \mbox{ \ \ and \ \ } E(u_0) <  (1-\delta_0)E(W).$$  
Let $I\ni 0$ be the maximal interval of existence given by Definition 2.10. Let $\overline{\delta} = \overline{\delta}(\delta_0,N)$ be as in Lemma 3.4.  Then, for each $t\in I$, we have
\begin{equation}\label{3.10}
\int\abso{\del u(t)}^2\le (1-\overline{\delta})\int\abso{\del W}^2
\end{equation}
\begin{equation}\label{3.11}
\int\abso{\del u(t)}^2 - \abso{u(t)}^{2^*}\ge \overline{\delta}\int\abso{\del u(t)}^2
\end{equation}
\begin{equation}\label{3.12}
E(u(t)) \ge 0.
\end{equation}
\end{theorem}
\begin{proof}
By Remark 2.9, $E(u(t)) = E(u_0)$, $t\in I$ and the Theorem follows directly from Lemma 3.4 and a continuity argument.
\end{proof}

\begin{corollary}\label{3.13}
Let $u,u_0$ be as in Theorem 3.9.  Then for all $t\in I$ we have $E(u(t))\simeq \int\abso{\del u(t)}^2\simeq \int\abso{\del u_0}^2$, with comparability constants which depend only on $\delta_0$.
\end{corollary}
\begin{proof}
$E(u(t))\le \int\abso{\del u(t)}^2$, but by (3.11) we have 
$$E(u(t))\ge \paren{\frac{1}{2} - \frac{1}{2^*}}\int\abso{\del u(t)}^2 + \frac{1}{2^*}\paren{\int \abso{\del u(t)}^2 - \abso{u(t)}^{2^*}}\ge C_{\overline{\delta}}\int\abso{\del u(t)}^2,$$
so the first equivalence follows.  For the second one note that $E(u(t)) = E(u_0)\simeq \int\abso{\del u_0}^2$, by the first equivalence when $t=0$.
\end{proof}

\begin{remark}\label{3.14}
Assume that $u_0\in \dot H^1$ and that 
$\abso{x}u_0\in L^2(\RR^N)$.  Assume that 
$$E(u_0) < E(W), \mbox{ \ \ but \ \ }\int\abso{\del u_0}^2 > \int\abso{\del W}^2. $$ 
If we choose $\delta_0$ so that $E(u_0)<(1-\delta_0) E(W)$, arguing as in Lemma 3.4 we can conclude that $\int\abso{\del u_0}^2 > (1+\overline{\delta})\int\abso{\del W}^2$, $\overline{\delta} = \overline{\delta}(\delta_0,N)$.  But then, 
$\int \abso{\del u_0}^2 - \abso{u_0}^{2^*} = 2^*E(u_0) - \paren{\frac{2}{N-2}}\int\abso{\del u_0}^2\le 2^*E(W) - \frac{2}{(N-2)} \frac{1}{C_N^N} - \frac{2\overline{\delta}}{(N-2)} \frac{1}{C_N^N} =  \frac{-2\overline{\delta}}{(N-2)C_N^N}$.
Since $E(u(t)) = E(u_0)$, a continuity argument shows that for all $t\in I$, the maximal interval of existence, we have $\int\abso{\del u(t)}^2\ge (1+\overline{\delta})\int\abso{\del W}^2$ and $\int\abso{\del u(t)}^2 - \int\abso{u(t)}^{2^*}\le \frac{-2\overline{\delta}}{(N-2)C_N^N}$.  But, the virial identity (\cite{G}) shows that, if $\abso{x}u_0\in L^2(\RR^N)$ then 
$\frac{d^2}{dt^2}\int \abso{x}^2\abso{u_0(x,t)}^2\,dx = 8\braces{\int\abso{\del u(t)}^2 - \abso{u(t)}^{2^*}}\le \frac{-16\overline{\delta}}{(N-2)C_N^N}$.  This shows that $I$ must be finite, i.e., the maximal interval of existence is finite. This argument is the critical analoge of the $H^1$ subcritical result in \cite{BC}.\\
Note that in the case where $u_0\in \dot H^1$ and 
$u_0\in L^2(\RR^N)$, the same result holds. Indeed,  one can use a local version of the virial identity (See section 5 for such a version) and the extra conservation law of the $L^2$ norm in time to control correction terms to obtain
$\frac{d^2}{dt^2}\int \phi(\abso{x})\abso{u_0(x,t)}^2\,dx \le \frac{-8\overline{\delta}}{(N-2)C_N^N}$, where $\phi$ is a regular and compactetly supported function (See for example Ogawa and Tsutsumi \cite{OT}).
\end{remark}

\section{Existence and compactness of a critical element}\label{4}
Let us consider the statement:
\begin{enumerate}[(SC)]
\item For all $u_0\in \dot H^1(\RR^N)$, with $\int\abso{\del u_0}^2 < \int\abso{\del W}^2$ and $E(u_0) < E(W)$, if $u$ is the corresponding solution to the (CP), with maximal interval of existence $I$ (see Definition 2.10), then $I = (-\infty, +\infty)$ and $\norm{S((-\infty,+\infty))}{u}<\infty$.
\end{enumerate}

We say that  $(SC)(u_0)$ holds if for this particular $u_0$, with $\int\abso{\del u_0}^2 < \int\abso{\del W}^2$ and $E(u_0) < E(W)$ and $u$ the corresponding solution to the (CP), with maximal interval of existence $I$ we have $I = (-\infty, +\infty)$ and $\norm{S((-\infty,+\infty))}{u}<\infty$.

Note that, because of Remark 2.7, if $\norm{L^2}{\del u_0}\le \overline{\delta}$, $(SC)(u_0)$ holds.  Thus, in light of Corollary 3.13, there exists $\eta_0>0$ such that, if $u_0$ is as in (SC) and $E(u_0)<\eta_0$, then $(SC)(u_0)$ holds.  Moreover, for any $u_0$ as in (SC), $E(u_0)\ge 0$, in light of Theorem 3.9.  Thus, there exists a number $E_C$, with $\eta_0\le E_C\le E(W)$, such that, if $u_0$ is as in (SC) and $E(u_0)<E_C$, $(SC)(u_0)$ holds and $E_C$ is optimal with this property.  For the rest of this section we will assume that $E_C<E(W)$.  
We now prove that there exits a critical element  $u_{0,C}$ at the critical level of energy $E_C $ so that $(SC)(u_{0,C})$ does not hold and from the minimality, this element has a compactness property up to the symetries of this equation. This is in fact a general principle which follows from the concentration compactness ideas. More precisely,

\begin{proposition}\label{4.1}
There exists $u_{0,C}$ in $\dot H^1$, with 
$$E(u_{0,C}) = E_C < E(W),  \int\abso{\del u_{0,C}}^2 < \int\abso{\del W}^2$$
  such that, if $u_C$ is the solution of (CP) with data $u_{0,C}$, and maximal interval of existence $I$, $0\in \mathring{I}$, then $\norm{S(I)}{u_C} = +\infty$.
\end{proposition}

\begin{proposition}\label{4.2}
Assume  $u_C$ is as in Proposition 4.1 and that (say) $\norm{S(I_+)}{u_C} = +\infty$, where $I_+ = (0,+\infty)\cap I$.  Then there exists $x(t)\in \RR^N$ and $\lambda(t)\in \RR^+$, for $t\in I_+$, such that 
$$K = \braces{v(x,t)~:~v(x,t) = \frac{1}{\lambda(t)^{(N-2)/2}}u_C\paren{\frac{x-x(t)}{\lambda(t)},t}}$$
has the property that $\overline{K}$ is compact in $\dot H^1$.  A corresponding conclusion is reached if $\norm{S(I_-)}{u_C} = +\infty$, where $I_- = (-\infty, 0)\cap I$.
\end{proposition}

The main tools that we will need in order to prove Propositions 4.1 and 4.2 are the following Lemmas.
\begin{lemma}[Concentration compactness]\label{4.3}
Let $\braces{v_{0,n}}\in \dot H^1$, $\int\abso{\del v_{0,n}}^2\le A$.  Assume that $\norm{L^{2(N+2)/N-2}}{e^{it\Delta}v_{0,n}}\ge \delta>0$, where $\delta = \delta(N)$ is as in Theorem 2.5.  Then there exists a sequence $\braces{V_{0,j}}_{j=1}^\infty$ in $\dot H^1$, a subsequence of $\braces{v_{0,n}}$ (which we still call $\braces{v_{0,n}}$) and a triple $(\lambda_{j,n};x_{j,n};t_{j,n})\in \RR^+\times\RR^N\times \RR$, with
\[
	\frac{\lambda_{j,n}}{\lambda_{j',n}} + \frac{\lambda_{j',n}}{\lambda_{j,n}} + \frac{\abso{t_{j,n}-t_{j',n}}}{\lambda_{j,n}^2} + \frac{\abso{x_{j,n}-x_{j',n}}}{\lambda_{j,n}}\to \infty
\]
as $n\to \infty$ for $j\ne j'$ (we say that $(\lambda_{j,n};x_{j,n};t_{j,n})$ is orthogonal if this property is verified) such that
\begin{equation}\label{4.4}
\norm{\dot H^1}{V_{0,1}}\ge \alpha_0(A)>0
\end{equation}
\begin{equation}\label{4.5}
\begin{aligned}
&\text{If $V_j^l(x,t) = e^{it\Delta}V_{0,j}$, then, given $\epsilon_0>0$, there exists $J = J(\epsilon_0)$ and} \\ &\text{$\braces{w_n}_{n=1}^\infty\in \dot H^1$, so that $v_{0,n} = \sum_{j=1}^J\frac{1}{\lambda_{j,n}^{(N-2)/2}}V_j^l\paren{\frac{x-x_{j,n}}{\lambda_{j,n}},\frac{-t_{j,n}}{\lambda_{j,n}^2}} + w_n$} \\ &\text{with $\norm{S((-\infty,+\infty))}{e^{it\Delta}w_n}\le \epsilon_0$, for $n$ large}
\end{aligned}
\end{equation}
\begin{equation}\label{4.6}
\int\abso{\del v_{0,n}}^2 = \sum_{j=1}^J\int\abso{\del V_{0,j}}^2 + \int\abso{\del w_n}^2 + o(1)\text{ as }n\to \infty
\end{equation}
\begin{equation}\label{4.7}
E(v_{0,n}) = \sum_{j=1}^JE(V_j^l(-t_{j,n}/\lambda_{j,n}^2))+E(w_n) + o(1)\text{ as }n\to \infty.
\end{equation}
\end{lemma}

\begin{remark}\label{4.8}
Lemma 4.3 is due to Keraani (\cite{K}).  It is based on the ``refined Sobolev inequality" ($N=3$)
\[
	\norm{L^6(\RR^3)}{h}\le C\norm{L^2(\RR^3)}{\del h}^{1/3}\norm{\dot B_{2,\infty}^0}{\del h}^{2/3},
\]
where $\dot B_{2,\infty}^0$ is the standard Besov space (\cite{BL}, \cite{GMO}).  (4.4) is a consequence of the proof of Corollary 1.9 in \cite{K}, (here, we use the hypothesis $\norm{L^{2(N+2)/N-2}}{e^{it\Delta}v_{0,n}}\ge\delta>0$) while (4.7) follows from the orthogonality of $(\lambda_{j,n};x_{j,n};t_{j,n})$ as in the proof of (4.6).  The rest of the Lemma is contained in the proof of Theorem 1.6 in \cite{K}. See also \cite{Ge}, \cite{MV}, \cite{BG}, \cite{K1}.
\end{remark}

\begin{lemma}\label{4.9}
Let $\braces{z_{0,n}}\in \dot H^1$, with $\int\abso{\del z_{0,n}}^2<\int\abso{\del W}^2$ and $E(z_{0,n})\to E_C$ and with $\norm{S((-\infty,+\infty))}{e^{it\Delta}z_{0,n}}\ge \delta$, with $\delta$ as in Theorem 2.5.  Let $\braces{V_{0,j}}$ be as in Lemma 4.3.  Assume that one of the two hypothesis
\begin{equation}\label{4.10}
\varliminf_{n\to \infty} E(V_1^l(-t_{1,n}/\lambda_{1,n}^2))<E_C
\end{equation}
or after passing to a subsequence, we have that, with $s_n = -t_{1,n}/\lambda_{1,n}^2$, $E(V_1^l(s_n))\to E_C$, and $s_n\to s_*\in [-\infty, +\infty]$, and if $U_1$ is the non-linear profile (see Definition 2.12 and Remark 2.13) associated to $(V_{0,1},\braces{s_n})$ we have that the maximal interval of existence of $U_1$ is $I = (-\infty, +\infty)$ and $\norm{S((-\infty,+\infty))}{U_1}<\infty$ and
\begin{equation}\label{4.11}
\varliminf_{n\to \infty}E(V_1^l(-t_{1,n}/\lambda_{1,n}^2))=E_C
\end{equation}
Then (after passing to a subsequence), for $n$ large, if $z_n$ is the solution of (CP) with data at $t=0$ equal to $z_{0,n}$, then  $(SC)(z_{0,n})$ holds.
\end{lemma}

Let us first assume the validity of Lemma 4.9 and use it (together with Lemma 4.3) to establish Propositions 4.1 and 4.2.

\begin{proof}[Proof of Proposition 4.1]
By the definition of $E_C$, and the assumption that $E_C<E(W)$, we can find $u_{0,n}\in \dot H^1$, with $\int\abso{\del u_{0,n}}^2<\int\abso{\del W}^2$, $E(u_{0,n})\to E_C$, and such that if $u_n$ is the solution of (CP) with data at $t=0$, $u_{0,n}$ and maximal interval of existence $I_n = (-T_-(u_{0,n}),T_+(u_{0,n}))$, then $\norm{S((-\infty,+\infty))}{e^{it\Delta}u_{0,n}}\ge \delta = \delta(N)>0$, where $\delta$ is as in Theorem 2.5 and $\norm{S(I_N)}{u_n} = +\infty$.  (Here we are also using Lemma 2.1 and Theorem 2.5.)  Note also that, since $E_C<E(W)$, there exists $\delta_0>0$ so that, for all $n$, we have $E(u_{0,n})\le (1-\delta_0)E(W)$.  Because of Theorem 3.9, we can find $\overline\delta$ so that $\int\abso{\del u_n(t)}^2\le (1-\overline\delta)\int\abso{\del W}^2$ for all $t\in I_n$, all $n$.  Apply now Lemma 4.3 for $\epsilon_0>0$ and Lemma 4.9.  We then have, for $J = J(\epsilon_0)$, that
\begin{equation}\label{4.12}
u_{0,n} = \sum_{j=1}^J \frac{1}{\lambda_{j,n}^{(N-2)/2}}V_j^l\paren{\frac{x-x_{j,n}}{\lambda_{j,n}},\frac{-t_{j,n}}{\lambda_{j,n}^2}} + w_n
\end{equation}
\begin{equation}\label{4.13}
\int\abso{\del u_{0,n}}^2 = \sum_{j=1}^J\int\abso{\del V_{0,j}}^2 + \int\abso{\del w_n}^2 + o(1)
\end{equation}
\begin{equation}\label{4.14}
E(u_{0,n}) = \sum_{j=1}^JE\paren{V_{j^l}\paren{\frac{-t_{j,n}}{\lambda_{j,n}^2}}} + E(w_n) + o(1).
\end{equation}
Note that because of (4.13) we have, for all $n$ large, that $\int\abso{\del w_n}^2\le (1-\overline\delta/2)\int\abso{\del W}^2$ and $\int\abso{\del V_{0,j}}^2\le (1-\overline\delta/2)\int\abso{\del W}^2$.  From Corollary 3.8 it now follows that \\$E(V_j^l(-t_{j,n}/\lambda_{jn}^2))\ge 0$ and $E(w_n)\ge 0$.  From this and (4.14) it  follows that $E(V_1^l(-t_{1,n}/\lambda_{1,n}^2))\le E(u_{0,n}) + o(1)$ and hence $\varliminf_{n\to \infty}E(V_1^l(-t_{1,n}/\lambda_{1,n}^2))\le E_C$.  If the left-hand side is strictly less than $E_C$, Lemma 4.9 gives us a contradiction with the choice of $u_{0,n}$, for $n$ large (after passing to a subsequence).  Hence, the left-hand side must equal $E_C$.  \\
Let then $U_1$ be the non-linear profile associated to $(V_1^l,\braces{s_n})$, with $s_n = -t_{1,n}/\lambda_{1,n}^2$ (after passing to a subsequence).  We first note that we must have $J=1$.  This is because (4.14) and $E(u_{0,n})\to E_C$, $E(V_1^l(-s_n))\to E_C$ now imply that $E(w_n)\to 0$ and $E(V_j^l(-t_{j,n}/\lambda_{j,n}^2)\to 0$, $j=2,...,J$.  Using (3.6) and the argument in the proof of Corollary 3.13, we have $C\sum_{j=2}^J\int\abso{\del V_j^l(-t_{j,n}/\lambda_{j,n}^2)}^2 + C\int\abso{\del w_n}^2\to 0$.  We then have, since $\int\abso{\del V_j^l(-t_{j,n}/\lambda_{j,n}^2)}^2 = \int\abso{\del V_{0,j}}^2$ that $V_{0,j} = 0$, $j=2,..,J$ and $\int\abso{\del w_n}^2\to 0$.  Hence (4.12) becomes $u_{0,n} = \frac{1}{\lambda_{1,n}^{(N-2)/2}}V_1^l\paren{\frac{x-x_{1,n}}{\lambda_{1,n}},s_n} + w_n$.  Let $v_{0,n} = \lambda_{1,n}^{(N-2)/2}u_{0,n}(\lambda_{1,n}(x+x_{1,n}))$ and note that scaling gives us that $v_{0,n}$ verifies the same hypothesis as $u_{0,n}$.  Moreover, $\widetilde w_n = \lambda_{1,n}^{(N-2)/2}w_n(\lambda_{1,n}(x+x_{1,n}))$ still verifies $\int\abso{\del \widetilde w_n}^2\to 0$.  Thus
\[
v_{0,n} = V_1^l(s_n) + \widetilde w_n,\quad \int\abso{\del \widetilde w_n}^2\to 0.
\]

Let us return to $U_1$, the non-linear profile associated to $(V_{0,1},\braces{s_n})$ and let $I_1 = (T_-(U_1),T_+(U_1))$ be its maximal interval of existence (see Remark 2.13).  Note that, by definition of non-linear profile, we have $\int\abso{\del U_1(s_n)}^2 = \int\abso{\del V_1^l(s_n)}^2+o(1)$ and $E(U_1(s_n)) = E(V_1^l(s_n))+o(1)$.  Note that in this case $E(V_1^l(s_n)) = E_C+o(1)$ and that $\int\abso{\del V_1^l(s_n)}^2 = \int\abso{\del V_{0,1}}^2 = \int\abso{\del u_{0,n}}^2 + o(1) < \int\abso{\del W}^2$ for $n$ large by Theorem 3.9.  Let's fix $\overline s\in I_1$.  Then $E(U_1(s_n)) = E(U_1(\overline s))$, so that 
$$E(U_1(\overline s)) = E_C.$$  
Moreover, $\int\abso{\del U_1(s_n)}^2<\int\abso{\del W}^2$ for $n$ large and hence by (3.10) $\int\abso{\del U_1(\overline s)}^2<\int\abso{\del W}^2$.  If $\norm{S(I_1)}{U_1} < +\infty$, Lemma 2.11 gives us that $I_1 = (-\infty, +\infty)$ and we then obtain a contradiction from Lemma 4.9.  Thus, 
$$\norm{S(I_1)}{U_1} = +\infty$$
and we then set $u_C = U_1$ (after a translation in time to make $\overline s=0$).
\end{proof}

\begin{proof}[Proof of Proposition 4.2]
We argue by contradiction.  For brevity of notation, let us set $u(x,t) = u_C(x,t)$.  If not, there exists $\eta_0>0$ and a sequence $\braces{t_n}_{n=1}^\infty$, $t_n\ge 0$ such that, for all $\lambda_0\in \RR^+$, $x_0\in \RR^N$, we have
\begin{equation}\label{4.15}
\norm{\dot H^1_x}{\frac{1}{\lambda_0^{(N-2)/2}}u\paren{\frac{x-x_0}{\lambda_0},t_n} - u(x,t_{n'})}\ge \eta_0,\quad \text{for }n\ne n'.
\end{equation}
Note that (after passing to a subsequence, so that $t_n\to \overline t\in [0,T_+(u_0)]$), we must have $\overline t = T_+(u_0)$, in view of the continuity of the flow in $\dot H^1$, as guaranteed by Theorem 2.5.  Note that, in view of Theorem 2.5 we must also have $\norm{S((0,+\infty))}{e^{it\Delta}u(t_n)}\ge \delta$. \\ 
Let us apply Lemma 4.3 to $v_{0,n} = u(t_n)$ with $\epsilon_0 > 0$.  We next prove that $J=1$.  In fact, if $\varliminf_{n\to \infty}E(V_1^l(-t_{1,n}/\lambda_{1,n}^2))<E_C$, since $\int\abso{\del u(t)}^2\le (1-\overline \delta)\int\abso{\del W}^2$ by Theorem 3.9, for all $t\in I_+$ and $E(u(t)) = E(u_0) = E_C<E(W)$, by Lemma 4.9 we obtain a contradiction.  Hence, we must have $\varliminf_{n\to \infty}E(V_1^l(-t_{1,n}/\lambda_{1,n}^2)) = E_C$.  The argument used in the proof of Proposition 4.1 now applies and gives $J=1$, $\int\abso{\del w_n}^2\to 0$.  Thus, we have
\begin{equation}\label{4.16}
u(t_n) = \lambda_{1,n}^{-(N-2)/2}V_1^l\paren{\frac{x-x_{1,n}}{\lambda_{1,n}},\frac{-t_{1,n}}{\lambda_{1,n}^2}} + w_n,\quad \int\abso{\del w_n}^2\to 0.
\end{equation}
Our next step is to show that $s_n = \frac{-t_{1,n}}{\lambda_{1,n}^2}$ must be bounded.  To see this note that
\[
	e^{it\Delta}u(t_n) = \lambda_{1,n}^{-(N-2)/2}V_1^l\paren{\frac{x-x_{1,n}}{\lambda_{1,n}},\frac{t-t_{1,n}}{(\lambda_{1,n})^2}} + e^{it\Delta}w_n.
\]
Assume that $t_{1,n}/\lambda_{1,n}^2\le -C_0$, $C_0$ a large positive constant.  Then, since \\$\norm{S((-\infty,+\infty))}{e^{it\Delta}w_n}< \delta/2$ for $n$ large, and
\[
\norm{S((0,+\infty))}{\lambda_{1,n}^{-(N-2)/2}V_1^l\paren{\frac{x-x_{1,n}}{\lambda_{1,n}},\frac{t-t_{1,n}}{(\lambda_{1,n})^2}}}\le \norm{S((C_0,+\infty))}{V_1^l(y,s)}\le \delta/2,
\]
for $C_0$ large, we get a contradiction.

If, on the other hand, $\frac{t_{1,n}}{\lambda_{1,n}^2}\ge C_0$, for a large positive constant $C_0$, $n$ large, we have
\[
\norm{S((-\infty,0))}{\lambda_{1,n}^{-(N-2)/2}V_1^l\paren{\frac{x-x_{1,n}}{\lambda_{1,n}},\frac{t-t_{1,n}}{(\lambda_{1,n})^2}}}\le \norm{S((-\infty,-C_0))}{V_1^l(y,s)}\le \delta/2,
\]
for $C_0$ large.  Hence, $\norm{S((-\infty,0))}{e^{it\Delta}u(t_n)}\le\delta$, for $n$ large and hence, Theorem 2.5 now gives $\norm{S((-\infty,t_n))}{u}\le 2\delta$, which, since $t_n\to T_+(u_0)$ gives us a contradiction.  Thus $\abso{t_{1,n}/\lambda_{1,n}^2}\le C_0$ and after passing to a subsequence, 
$$t_{1,n}/\lambda_{1,n}^2\to t_0\in (-\infty,+\infty).$$  
But then, since (4.15) and (4.16) imply that, for $n \ne n'$ large (independently of $\lambda_0,x_0$) we have
\[\begin{aligned}
&\left|\kern-.1em\left|\frac{1}{(\lambda_0)^{(N-2)/2}}\frac{1}{(\lambda_{1,n})^{(N-2)/2}}V_1^l\paren{\frac{\frac{x-x_0}{\lambda_0}-x_{1,n}}{\lambda_{1,n}},-t_{1,n}/(\lambda_{1,n})^2}\right.\right. \\
&\left.\left.- \frac{1}{(\lambda_{1,n'})^{(N-2)/2}}V_1^l\paren{\frac{x-x_{1,n'}}{\lambda_{1,n'}},-t_{1,n'}/(\lambda_{1,n'})^2}\right|\kern-.1em\right|_{\dot H^1}\ge \eta_0/2
\end{aligned}\]
or
\[\begin{aligned}
&\left|\kern-.1em\left|
\paren{\frac{\lambda_{1,n'}}{\lambda_{1,n}\lambda_0}}^{(N-2)/2}
V_1^l\paren{\frac{y\lambda_{1,n'}}{\lambda_0\lambda_{1,n}}+\widetilde x_{n,n'}-\widetilde x_0,-t_{1,n}/(\lambda_{1,n})^2} \right.\right. \\
&-V_1^l\paren{y,-t_{1,n'}/\lambda_{1,n'}^2}\Big|\kern-.1em\Big|_{\dot H^1}\ge \eta_0/2,
\end{aligned}\]
where $\widetilde x_{n,n'}$ is a suitable point in $\RR^N$ and $\lambda_0,\widetilde x_0$ are arbitrary.  But if we choose $\lambda_0 = \lambda_{1,n'}/\lambda_{1,n}$, $\widetilde x_0 = x_{n,n'}$, we reach a contradiction since $-t_{1,n}/(\lambda_{1,n})^2\to -t_0$ and $-t_{1,n'}/(\lambda_{1,n'})^2\to -t_0$.
\end{proof}

Thus, to complete the proofs of Propositions 4.1 and 4.2 we only need to provide the proof of Lemma 4.9.
\begin{proof}[Proof of Lemma 4.9]
Let us assume first that (4.11) holds and set $A = \int\abso{\del W}^2$, $A' = \int\abso{\del W}^2$, $M = \norm{S((-\infty,+\infty))}{U_1}$.  Arguing (for some $\epsilon_0>0$ in Lemma 4.3) as in the proof of Proposition 4.1, we see that $\varliminf_{n\to \infty}E(V_1^l(-t_{1,n}/\lambda_{1,n}^2)) = E_C$ and $E_C<E(W)$, imply that $J=1$, $\int\abso{\del w_n}^2\to 0$.  Moreover, if 
$$v_{0,n} = \lambda_{1,n}^{(N-2)/2}z_{0,n}(\lambda_{1,n}(x+x_{1,n})), \widetilde w_n = \lambda_{1,n}^{(N-2)/2}w_n(\lambda_{1,n}(x+x_{1,n})), s_n = -\frac{t_{1,n}}{\lambda_{1,n}^2},$$ 
we have $\int\abso{\del \widetilde w_n}^2\to 0$ and $v_{0,n} = V_1^l(s_n)+\widetilde w_n$, while $\norm{S((-\infty,+\infty))}{e^{it\Delta}v_{0,n}}\ge \delta$, $\int\abso{\del v_{0,n}}^2 < \int\abso{\del W}^2$, $E(v_{0,n})\to E_C$.  Note now that $\int\abso{\del \widetilde V_1^l(s_n) - U_1(s_n)}^2 = o(1)$ by definition of non-linear profile.  We then have 
$$v_{0,n} = U_1(s_n) + \widetilde{\widetilde{w}}_n, \ \ \int\abso{\del \widetilde{\widetilde{w}}_n}^2\to 0.$$  
Moreover, as in the proof of Proposition 4.1, $E(U_1(0)) = E_C$ and $\int\abso{\del U_1(t)}^2<\int\abso{\del W}^2$ for all $t$.  We now apply Theorem 2.14, with $\epsilon_0<\epsilon_0(M,A,A',N)$ and $n$ large, with $\widetilde u = U_1$, $e\equiv 0$, $t_0 = 0$, $u_0 = v_{0,n}$.  This case now follows.

Assume next that (4.10) holds.  The first claim is that for $j\ge 2$ we also have $\varliminf_{n\to \infty}E(V_j^l(-t_{j,n}/\lambda_{j,n}^2))<E_C$.  In fact, after passing to a subsequence, assume $\lim_{n\to \infty}E(V_1^l(-t_{1,n}/\lambda_{1,n}))<E_C$.  Because of (4.6) we have
\[
	\int\abso{\del z_{0,n}}^2\ge \sum_{j=1}^J\int\abso{\del V_{0,j}}^2+o(1),
\]
and since $E_C < E(W)$, for $n$ large we have $E(z_{0,n})\le (1-\delta_0)E(W)$, by Lemma 3.4, $\int\abso{\del z_{0,n}}^2\le (1-\overline \delta)\int\abso{\del W}^2$ and hence $\int\abso{\del V_{0,j}}^2\le (1-\overline\delta)\int\abso{\del W}^2$.  Similarly, $\int\abso{\del w_n}^2\le (1-\overline\delta)\int\abso{\del W}^2$.  By Corollory 3.8, we have $E(V_j^l(-t_{j,n}/\lambda_{j,n}^2))\ge 0$, $E(w_n)\ge 0$.  Also, from (4.4) and the proof of Corollary 3.13, we have, for $n$ large, that $E(V_1^l(-t_{1,n}/\lambda_{1,n}^2))\ge C\int\abso{\del V_{0,1}}^2\ge c\alpha_0 = \overline{\alpha_0}>0$, so that, from (4.7) we obtain, for $n$ large
\[
	E(z_{0,n})\ge \overline{\alpha_0} + \sum_{j=2}^JE(V_j^l(-t_{j,n}/\lambda_{j,n}^2))+o(1),
\]
so that the claim follows from $E(z_{0,n})\to E_C$.  \\
We next claim that (after passing to a subsequence so that, for each $j$, \\
$\lim_nE(V_j^l(-t_{j,n}/\lambda_{j,n}^2))$ exists and $\lim_n(-t_{j,n}/\lambda_{j,n}^2) = \overline{s_j}\in [-\infty,+\infty]$ exists) if $U_j$ is the non-linear profile associated to $(V_j^l,\braces{-t_{j,n}/\lambda_{j,n}^2})$, then $U_j$ verifies (SC).  In fact, by definition of non-linear profile, $E(U_j) <E_C$, since $\lim_n E(V_j^l(-t_{j,n}/\lambda_{j,n}^2))<E_C$.  Moreover, since $\int\abso{\del V_j^l(-t_{j,n}/\lambda_{j,n}^2)}^2\le (1-\overline\delta)\int\abso{\del W}^2$, by the definition of non-linear profile and Theorem 3.9, if $\overline t\in I_j$, the maximal interval for $U_j$, $\int\abso{\del U_j(\overline t)}^2<\int\abso{\del W}^2$ so that, by the definition of $E_C$ our claim follows.  Note that the argument in the proof of Theorem 2.14 also gives that $\norm{W((-\infty,+\infty))}{\del U_j} < +\infty$.  \\
Our final claim is that there exists $j_0$ so that, for $j\ge j_0$ we have
\begin{equation}\label{4.17}
\norm{S((-\infty,+\infty))}{U_j}^{2(N+2)/N-2}\le C(\int\abso{\del V_{0,j}}^2)^{N+2/N-2}.
\end{equation}
In fact, from (4.6), for fixed $J$ we see that (choosing $n$ large) $\sum_{j=1}^J \int\abso{\del V_{0,j}}^2\le \int\abso{\del z_{0,n}}^2+o(1)\le 2\int\abso{\del W}^2$.  Thus, for $j\ge j_0$, we have $\int\abso{\del V_{0,j}}^2\le\widetilde \delta$, where $\widetilde \delta$ is so small that $\norm{S((-\infty,+\infty))}{e^{it\Delta}V_{0,j}}\le \delta$, with $\delta$ as in Theorem 2.5.  From Remark 2.13 it then follows that $\norm{S((-\infty,+\infty))}{U_j}\le 2\delta$, and using the integral equation in Remark 2.13, that $\norm{\dot H^1}{U_j(0)}\le C\norm{\dot H^1}{V_{0,j}}$ and $\norm{W((-\infty,+\infty))}{\del U_j}\le C\norm{\dot H^1}{V_{0,j}}$, which gives (4.17).  \\
For $\epsilon_0>0$, to be chosen, define now
\begin{equation}\label{4.18}
H_{n,\epsilon_0} = \sum_{j=1}^{J(\epsilon_0)}\frac{1}{\lambda_{j,n}^{(N-2)/2}}U_j\paren{\frac{x-x_{j,n}}{\lambda_{j,n}},\frac{t-t_{j,n}}{\lambda_{j,n}^2}}.
\end{equation}
We then have:
\begin{equation}\label{4.19}
\norm{S((-\infty,+\infty))}{H_{n,\epsilon_0}}\le C_0,
\end{equation}
uniformly in $\epsilon_0$, for $n\ge n(\epsilon_0)$.  In fact,
\[\begin{aligned}
\norm{S((-\infty,+\infty))}{H_{n,\epsilon_0}}^{2(N+2)/N-2}&
= \int\int \left[\sum_{j=1}^{J(\epsilon_0)}\frac{1}{\lambda_{j,n}^{(N-2)/2}}U_j\paren{\frac{x-x_{j,n}}{\lambda_{j,n}},\frac{t-t_{j,n}}{\lambda_{j,n}^2}}\right]^{\frac{2(N+2)}{N-2}} \\
&\le \sum_{j=1}^{J(\epsilon_0)}\int\int \abso{\frac{1}{\lambda_{j,n}^{(N-2)/2}}U_j\paren{\frac{x-x_{j,n}}{\lambda_{j,n}},\frac{t-t_{j,n}}{\lambda_{j,n}^2}}}^{\frac{2(N+2)}{N-2}} \\
&+ C_{J(\epsilon_0)}\sum_{j'\ne j}\int\int \abso{\frac{1}{\lambda_{j,n}^{(N-2)/2}}U_j\paren{\frac{x-x_{j,n}}{\lambda_{j,n}},\frac{t-t_{j,n}}{\lambda_{j,n}^2}}} \\
&\abso{\frac{1}{\lambda_{j',n}^{(N-2)/2}}U_{j'}\paren{\frac{x-x_{j',n}}{\lambda_{j',n}},\frac{t-t_{j',n}}{\lambda_{j',n}^2}}}^{\frac{N+6}{N-2}} 
= \mathrm{I+II}.
\end{aligned}\]
For $n$ large, $\mathrm{II}\to 0$, by the orthogonality of $(\lambda_{j,n};x_{j,n};t_{j,n})$ (see Keraani \cite{K}, Lemma 2.7, (2.95), (2.96), etc.)  Hence, for $n$ large we have $\mathrm{II\le I}$.  But (with $j_0$ as in (4.17)),
\[\begin{aligned}
I
&\le \sum_{j=1}^{j_0}\norm{S((-\infty,+\infty))}{U_j}^{2(N+2)/N-2} + \sum_{j=j_0}^{J(\epsilon_0)}\norm{S((-\infty,+\infty))}{U_j}^{2(N+2)/N-2}  \\
&\le \sum_{j=1}^{j_0} \norm{S((-\infty,+\infty))}{U_j}^{2(N+2)/N-2} + C\sum_{j=j_0}^{J(\epsilon_0)}(\int\abso{\del V_{0,j}}^2)^{N+2/N-2}\le C_0/2
\end{aligned}\]
because of (4.6).

For $\epsilon_0>0$, to be chosen, define
\begin{equation}\label{4.20}
\begin{aligned}
&R_{n,\epsilon_0}
= \abso{H_{n,\epsilon_0}}^{\frac{4}{N-2}} H_{n,\epsilon_0} \\
&- \sum_{j=1}^{J(\epsilon_0)}\abso{\frac{1}{\lambda_{j,n}^{(N-2)/2}}U_j\paren{\frac{x-x_{j,n}}{\lambda_{j,n}},\frac{t-t_{j,n}}{\lambda_{j,n}^2}}}^{\frac{4}{N-2}} 
\frac{1}{\lambda_{j,n}^{(N-2)/2}}U_j\paren{\frac{x-x_{j,n}}{\lambda_{j,n}},\frac{t-t_{j,n}}{\lambda_{j,n}^2}}.
\end{aligned}
\end{equation}
We then have
\begin{equation}\label{4.21}
\text{For $n = n(\epsilon_0)$ large, \ \ } \norm{L_t^2L_x^{2N/N+2}}{\del R_{n,\epsilon_0}}\to 0\quad\text{as $n\to \infty$}.
\end{equation}
This follows from the orthogonality of $(\lambda_{j,n};x_{j,n};t_{j,n})$, the fact that $\norm{S((-\infty,+\infty))}{U_j}<\infty$, $\norm{W((-\infty,+\infty))}{\del U_j}<\infty$, and arguments of Keraani \cite{K} (see in particular (2.95), (2.96)).

We now will apply Theorem 2.14.  Let $\widetilde u = H_{n,\epsilon_0}$, $e = R_{n,\epsilon_0}$, where $\epsilon_0$ is still to be determined.  Recall that $z_{0,n} = \sum_{j=1}^{J(\epsilon_0)}\frac{1}{\lambda_{j,n}^{(N-2)/2}}V_j^l\paren{\frac{x-x_{j,n}}{\lambda_{j,n}},\frac{-t_{j,n}}{\lambda_{j,n}^2}} + w_n$, where $\norm{S((-\infty,+\infty))}{e^{it\Delta}w_n}\le \epsilon_0$.  By the definition of non-linear profile, we now have
\begin{equation}\label{4.22}
z_{0,n}(x) = H_{n,\epsilon_0}(x,0) + \widetilde w_n(x),
\end{equation}
where, for $n$ large $\norm{S((-\infty,+\infty))}{e^{it\Delta}\widetilde w_n}\le 2\epsilon_0$.  \\
Notice also that, because of the orthogonality of $(\lambda_{j,n};x_{j,n};t_{j,n})$, for $n = n(\epsilon_0)$ large, we have (using also Corollary 3.13), that \\
$\int\abso{\del H_{n,\epsilon_0}(t)}^2\le 2\sum_{j=1}^{J(\epsilon_0)}E_0\paren{U_j\paren{\frac{t-t_{j,n}}{\lambda_{j,n}^2}}}\le 4C\sum_{j=1}^{J(\epsilon_0)}\int\abso{\del V_{0,j}}^2$, and \\
$\sum_{j=1}^{J(\epsilon_0)}\int\abso{\del V_{0,j}}^2\le \int\abso{\del z_{0,n}}^2\int\abso{\del z_{0,n}}^2+o(1)\le 2\int\abso{\del W}^2$.  \\
Let now $M = C_0$, with $C_0$ as in (4.19), $A = \widetilde{C} \int\abso{\del W}^2$, $A' = A + \int\abso{\del W}^2$, $\epsilon_0 < \epsilon_0(M,A,A',N)/2,$ where $\epsilon_0(M,A,A',N)$ is as in Theorem 2.14.  Fix $\epsilon_0$ and choose $n$ so large that $\norm{L_t^2L_x^{2N/N+2}}{\del R_{n,\epsilon_0}}<\epsilon_0$ and so that all the above properties hold.  Then Theorem 2.14 gives the conclusion of Lemma 4.9 in the case when (4.10) holds.
\end{proof}

\begin{remark}\label{4.23}
Assume that $\braces{z_{0,n}}$ in Lemma 4.3 are all radial.  Then $V_{0,j}$, $w_n$ can be chosen to be radial and we can choose $x_{j,n}\equiv 0$.  This follows directly from Keraani's proof \cite{K}.  If we then define (SC) and $E_C$ by restricting only to radial functions, we obtain a $u_C$ as in Proposition 4.1 which is radial, and we can establish Proposition 4.2 with $x(t)\equiv 0$.
\end{remark}

\section{Rigidity Theorem}\label{5}
In this section we will prove the following:
\begin{theorem}\label{5.1}
Assume that $u_0\in \dot H^1$ is such that 
$$E(u_0)<E(W), \ \ \int\abso{\del u_0}^2<\int\abso{\del W}^2.$$  
Let $u$ be the solution of (CP) with $u|_{t=0} = u_0$, with maximal interval of existence $(-T_-(u_0),T_+(u_0))$ (see Definition 2.10). Assume that there exists $\lambda(t)>0$, for $t\in [0,T_+(u_0))$, with the property that
\[
K = \braces{v(x,t) = \frac{1}{\lambda(t)^{(N-2)/2}}u\paren{\frac{x}{\lambda(t)},t}~:~t\in [0,T_+(u_0))}
\]
is such  that $\overline K$ is compact in $\dot H^1$.  Then $T_+(u_0) = +\infty$, $u_0 \equiv  0$.
\end{theorem}

\begin{remark}\label{5.2}
We conjecture that Theorem 5.1 remains true if \\ $v(x,t) = \frac{1}{\lambda(t)^{(N-2)/2}}u\paren{\frac{x-x(t)}{\lambda(t)},t}$, with $x(t)\in \RR^N$, $t\in [0,T_+(u_0))$.   In other words, for ``energy subcritical'' initial data, compactness up to the invariances of the equation, for solutions, is only true for $u \equiv  0$.
\end{remark}
We start out with a special case of the strengthened form of Theorem 5.1, namely:

\begin{proposition}\label{5.3}
Assume that $u$, $v$, $\lambda(t)$, $x(t)$ are as in Remark 5.2, that $\abso{x(t)}\le C_0$ and that $\lambda(t)\ge A_0>0$.  Then the conclusion of Theorem 5.1 holds.  Moreover, if $T_+(u_0)<+\infty$, the hypothesis $\abso{x(t)}\le C_0$ is not needed.
\end{proposition}

\begin{remark}\label{5.4}
Because of the continuity of $u(t)$ in $\dot H^1$, it is clear that in proving Proposition 5.3 we can assume that $\lambda(t),x(t)\in C^\infty([0,T_+(u_0)))$ and that $\lambda(t)>0$ for each $t\ge 0$. Indeed, first  by the compactness of $\overline K$ and the theory of (CP), we construct  piecewise contant  (with small jumps) $\lambda_1(t),x_1(t)$ such that the corresponding set $K_1$  is included in  $\tilde{K}_1 = \braces{w(t) \mbox{ \ solution of (CP) with initial data in \ }  \overline K, t\in[0,t_0]}$, $t_0$ small. 
Then we can contruct regular  $\lambda_2(t),x_2(t)$ such that $K_2$ is included in the precompact set
$\braces{ \lambda_0^{-\frac{(N-2)}{2}}w((x-x_0)\lambda_0^{-1}),  \mbox{ \ for \ }  w \in \tilde{K}_1, 1/2 \le \lambda_0 \le 2, |x_0|\le1}.$ 
 The continuity of $\lambda(t),x(t)$ will not be used in our proof.
\end{remark}

In the next lemma we will collect some useful facts:
\begin{lemma}\label{5.5}
Let $u,v$ be as in Proposition 5.3.  
\begin{enumerate}[i)]
\item Let $\delta_0>0$ be such that $E(u_0)\le (1-\delta_0)E(W)$.  Then for all $t\in [0,T_+(u_0))$, we have
\begin{gather*}
\int\abso{\del u(t)}^2\le (1-\overline \delta)\int\abso{\del W}^2; \\
\shoveleft \int \abso{\del u(t)}^2-\abso{u(t)}^{2^*}\ge \overline \delta\int \abso{\del u(t)}^2; \\
\shoveleft C_{1,\delta_0}\int\abso{\del u_0}^2\le E(u_0)\le C_{2}\int\abso{\del u_0}^2, \\
E(u(t)) = E(u_0), \\
C_{1,\delta_0}\int\abso{\del u_0}^2\le \int\abso{\del u(t)}^2\le C_{2}\int\abso{\del u_0}^2.
\end{gather*}
\item
\[
	\int\abso{\del v(t)}^2\le C_{2}\int\abso{\del W}^2
\]
and
\[
	\norm{L_x^{2^*}}{v(t)}^2\le C_{3}\int\abso{\del W}^2.
\]
\item For all $x_0\in \RR^N$
\[
	\int\frac{\abso{v(x,t)}^2}{\abso{x-x_0}^2}\le C_{4}\int\abso{\del W}^2.
\]
\item For each $\epsilon>0$, there exists $R(\epsilon_0)>0$, such that, for $0\le t < T_+(u_0)$, we have
\[
	\int_{|x|\ge R(\epsilon_0)} \abso{\del v}^2 + \abso{v}^{2^*} + \frac{\abso{v}^2}{\abso{x}^2}\le \epsilon_0.
\]
\end{enumerate}
\end{lemma}
\begin{proof}
i) follows from Theorem 3.9 and Corollary 3.13.  ii) follows from i) by Sobolev embedding, while iii) follows from i) by Hardy's inequality.  iv) follows (using Sobolev embedding and the Hardy inequality) from the compactness of $\overline K$.
\end{proof}

The next lemma is a localized virial identity, in the spirit of Merle \cite{M}, Lemma 3.6.
\begin{lemma}\label{5.6}
Let $\varphi\in C_0^\infty(\RR^N)$, $t\in [0,T_+(u_0))$.  Then:
\begin{enumerate}[i)]
\item
\[
\frac{d}{dt}\int\abso{u}^2\varphi\,dx = 2\im\int\overline u\del u\del \varphi\,dx
\]
\item
\[\begin{aligned}
\frac{d^2}{dt^2}\int\abso{u}^2\varphi\,dx
&= 4\sum_{l,j}\re\int\partial_{x_l}\partial_{x_j}\varphi\cdot\partial_{x_l}u\cdot \partial_{x_j}\overline u  \\
&- \int \Delta^2\varphi\abso{u}^2 - \frac{4}{N}\int \Delta\varphi\abso{u}^{2^*}. 
\end{aligned}\]
\end{enumerate}
\end{lemma}
The proof of Lemma 5.6 is standard, see \cite{M} and Glassey \cite{G}.

\begin{proof}[Proof of Proposition 5.3] The proof splits in two cases, the finite time blow-up case for $u$ and the infinite time of existence for $u$.

\underline{Case 1}:  $T_+(u_0)<\infty$.  (In this case we don't need the assumption $\abso{x(t)}<C_0$ or the energy constraints on $u$, only $sup_{t \in [0,T_+(u_0))} \int\abso{\del u(t)}^2 < \infty$ is needed. Note that this rules out the existence of self similar solutions in $ \dot H^1$, i.e. solutions for which $\lambda(t)\sim (T-t)^{-1/2}$.)  \\
Note first that $\lambda(t)\to \infty$ as $t\to T_+(u_0)$.  If not, there exists $t_i\uparrow T_+(u_0)$, with $\lambda(t_i)\to \lambda_0<+\infty$.  Let $v_i(x) = \frac{1}{\lambda(t_i)^{(N-2)/2}}u\paren{\frac{x-x(t_i)}{\lambda(t_i)},t_i}$ and let $v(x)\in \dot H^1$ be such that $v_i\to v$ in $\dot H^1$ (from the compactness of $\overline K$).  Hence, $u\paren{x-\frac{x(t_i)}{\lambda(t_i)},t_i} = \lambda(t_i)^{(N-2)/2}v_i(\lambda(t_i)x)\to \lambda_0^{(N-2)/2}v(\lambda_0x)$ in $\dot H^1$ (since $\lambda(t_i)\ge A_0$, $\lambda_0\ge A_0$).  Let now $h(x,t)$ be the solution of (CP), given by Remark 2.8 with data $\lambda_0^{(N-2)/2}v(\lambda_0x)$ at time $T_+(u_0)$, in an interval $(T_+(u_0)-\delta,T_+(u_0)+\delta)$, wtih $\norm{S((T_+(u_0)-\delta,T_+(u_0)+\delta))}{h}<\infty$.  Let $h_i(x,t)$ be the solution with data at $T_+(u_0)$ equal to $u\paren{x-\frac{x(t_i)}{\lambda(t_i)},t_i}$.  Then, the (CP) theory guarantees that
\[
	\sup_i\norm{S((T_+(u_0)-\frac{\delta}{2},T_+(u_0)+\frac{\delta}{2}))}{h_i}<\infty.
\]
But, $u(x- \frac{x(t_i)}{\lambda(t_i)},t+t_i-T_+(u_0)) = h_i(x,t)$, contradicting Lemma 2.11, since $T_+(u_0)<\infty$.  \\
Let us prove now a decay result for $u$ from the concentration properties in $L^{2^*}$ of $u$ at $T_+(u_0)$. Let us now fix $\varphi\in C_0^\infty(\RR^N)$, $\varphi$ radial, $\varphi\equiv 1$ for $\abso{x}\le 1$, $\varphi\equiv 0$ for $\abso{x}\ge 2$ and set $\varphi_R(x) = \varphi(x/R)$.  Define
\begin{equation}\label{5.7}
y_R(t) = \int \abso{u(x,t)}^2\varphi_R(x)\,dx,\quad t\in [0,T_+(u_0)).
\end{equation}
We then have:
\begin{equation}\label{5.8}
\abso{y_R'}\le C_N \int\abso{\del W}^2.
\end{equation}
In fact, by Lemma 4.6, i)
\[\begin{aligned}
\abso{y_R'}
&\le \frac{2}{R}\abso{\im \int\overline u\del u\del \varphi(x/R)\,dx} \\
&\le C_N\paren{\int \abso{\del u}^2}^{1/2}\cdot \paren{\int\frac{\abso{u}^2}{\abso{x}^2}}^{1/2}\le C_N\int\abso{\del W}^2,
\end{aligned}\]
by ii) in Lemma 5.5.

We also have:
\begin{equation}\label{5.9}
\text{For all $R>0$,}\quad \int_{|x|<R}\int\abso{u(x,t)}^2\,dx\to 0\quad\text{as }t\to T_+(u_0).
\end{equation}
In fact, $u(y,t) = \lambda(t)^{(N-2)/2}v(\lambda(t)y+x(t),t)$, so that
\[\begin{aligned}
\int_{|x|<R}\abso{u(x,t)}^2\,dx
&= \lambda(t)^{-2}\int_{|y|<R\lambda(t)}\abso{v(y+x(t),t)}^2\,dy \\
&= \lambda(t)^{-2}\int_{B(x(t),R\lambda(t))} \abso{v(z,t)}^2\,dz \\
&= \lambda(t)^{-2}\int_{B(x(t),R\lambda(t))\cap B(0,\epsilon R\lambda(t))}\abso{v(z,t)}^2\,dz \\
&+ \lambda(t)^{-2}\int_{B(x(t),R\lambda(t))\setminus B(0,\epsilon R\lambda(t))}\abso{v(z,t)}^2\,dz \\
&=A+B.
\end{aligned}\]
By H\"older's inequality, $A\le \lambda(t)^{-2}(\epsilon R\lambda(t))^{N \frac2N} \norm{L^{2^*}}{v}^2 \le \epsilon^2R^2C_{3}\int\abso{\del W}^2$, which is small with $\epsilon$.
\[
	B\le \lambda(t)^{-2}(R\lambda(t))^{N \frac2N}\norm{L^{2^*}{(\abso{x}\ge \epsilon R\lambda(t))}}{v}^2 
= R^2\norm{L^{2^*}{(\abso{x}\ge \epsilon R\lambda(t))}}{v}^2 \to 0\quad\text{as }t\to T_+(u_0),
\]
by iv) in Lemma 5.5, since $\lambda(t)\to +\infty$ as $t\to T_+(u_0)$.  \\
From (5.9) and (5.8), we have:
$$
y_R(0)
\le y_R(T_+(u_0)) + C_NT_+(u_0) \int\abso{\del W}^2 = C_NT_+(u_0)\int\abso{\del W}^2.
$$
Thus, letting $R \to +\infty$ we obtain $$u_0\in L^2(\RR^N).$$  
Arguing as before, $\abso{y_R(t) - y_R(T_+(u_0))}\le C_N(T_+(u_0)-t)\int\abso{\del W}^2$ so that $y_R(t)\le C_N(T_+(u_0) - t)\int\abso{\del W}^2$.  Letting $R\to \infty$, we see that $\norm{L^2}{u(t)}^2\le C_N(T_+(u_0)-t)\int\abso{\del W}^2$ and so by the conservation of the $L^2$ norm $\norm{L^2}{u(T_+(u_0))} = \norm{L^2}{u_0} = 0$.  But then $u\equiv 0$, contradicting $T_+(u_0)<\infty$.

\underline{Case 2}: $T_+(u_0) = +\infty$. \\
 In this case we assume, in addition, that $\abso{x(t)}\le C_0$.  We first note that
\begin{gather}\label{5.10}
\text{For each $\epsilon>0$, there exists $R(\epsilon)>0$ such that, for all $t\in [0,\infty)$, we have:} \\
\int_{|x|>R(\epsilon)}\abso{\del u}^2 + \abso{u}^{2^*} + \frac{\abso{u}^2}{\abso{x}^2}\le \epsilon. \notag
\end{gather}
In fact, $u(y,t) = \lambda(t)^{(N-2)/2}v(\lambda(t)y+x(t),t)$, so that
\[\begin{aligned}
\int_{|y|>R(\epsilon)}\abso{\del u(y,t)}^2\,dy
&= \int_{|y|>R(\epsilon)}\lambda(t)^N\abso{\del v(\lambda(t)y+x(t),t)}^2\,dy\\
= \int_{|z|>R(\epsilon)\lambda(t)}\abso{\del v(z+x(t),t)}^2\,dz
&\le \int_{|z|\ge R(\epsilon)A_0}\abso{\del v(z+x(t),t)}^2\,dz  \\
&\le \int_{|\alpha|\ge R(\epsilon)A_0-C_0}\abso{\del v(\alpha,t)}^2\,d\alpha,
\end{aligned}\]
and the statement for this term now follows from Lemma 5.5 iv).  The other terms are handled similarly.
\begin{gather}\label{5.11}
\text{There exists $R_0>0$ such that, for all $t\in [0,+\infty)$, we have} \\
8\int_{|x|\le R_0}\abso{\del u}^2 - 8\int_{|x|\le R_0}\abso{u}^{2^*}\ge C_{\delta_0}\int\abso{\del u_0}^2 \notag.
\end{gather}
In fact, (3.11) combined with Lemma 5.5 i) yields $8\int \abso{\del u}^2 - 8\int \abso{u}^{2^*}\ge \widetilde C_{\delta_0}\int\abso{\del u_0}^2$.  Now combine this with (5.10), with $\epsilon = \epsilon_0\int\abso{\del u_0}^2$ to obtain (5.11).  \\
To prove Case 2, we choose $\varphi\in C_0^\infty(\RR^N)$, radial, with $\varphi(x) = \abso{x}^2$ for $\abso{x}\le 1$, $\varphi(x)\equiv 0$ for $\abso{x}\ge 2$.  Define $z_R(t) = \int\abso{u(x,t)}^2R^2\varphi\paren{\frac{x}{R}}\,dx$.  We then have:
$$ 
\mbox{for \ } t>0, \ \ \abso{z_R'(t)}\le C_{N,\delta_0}\int\abso{\del u_0}^2R^2
$$
$$
\mbox{for $R$ large enough, $t>0,$ \ \ } 
 z_R''\ge C_{N,\delta_0}\int\abso{\del u_0}^2.
$$
In fact, from Lemma 5.6, i),
\[\begin{aligned}
\abso{z_R'(t)}
&\le 2R\abso{\im \int\overline u \del u\del \varphi\paren{\frac{x}{R}}\,dx} \le  
C_NR \int_{0\le|x|\le2R}\frac {|x|} {|x|}\abso{\del u}\abso{u}\\
&\le C_NR^2 \paren{\int\abso{\del u}^2}^{1/2}\paren{\int\frac{\abso{u}^2}{\abso{x}^2}}^{1/2}\le C_{N}R^2\int\abso{\del u_0}^2,
\end{aligned}\]
because of Lemma 5.5 i), while from Lemma 5.6, ii),
\begin{multline*}
z_R''(t)
= 4\sum_{l,j}\re \int\partial_{x_l}\partial_{x_j}\varphi\paren{\frac{x}{R}}\partial_{x_l}u\cdot \partial_{x_j}\overline u - \int\Delta^2\varphi\paren{\frac{x}{R}}\frac{\abso{u}^2}{R^2}  \\
- \frac{4}{N}\int \Delta\varphi\abso{u}^{2^*}\ge 8\left[\int_{|x|\le R}\abso{\del u}^2-\abso{u}^{2^*}\right] \\
-C_N\int_{R\le |x|\le2R}\left[ \abso{\del u}^2 + \frac{\abso{u}^2}{\abso{x}^2} + \abso{u}^{2^*}\right]\ge C_{N,\delta_0}\int\abso{\del u_0}^2
\end{multline*}
for $R$ large, in view of (5.11) and (5.10).

If we now integrate in $t$, we have $z_R'(t) - z_R'(0)\ge C_{N,\delta_0}t\int\abso{\del u_0}^2$, but we also have $\abso{z_R'(t)-z_R'(0)}\le 2C_{N}R^2\int\abso{\del u_0}^2$, a contradiction for $t$ large, unless $\int\abso{\del u_0}^2 = 0$.
\end{proof}

\begin{proof}[Proof of Theorem 5.1]
(See \cite{M2} for similar proof) Assume that $u_0\not\equiv 0$ so that $\int\abso{\del u_0}^2>0$
and because of Lemma 5.5 i) (which is still valid here), $E(u_0)\ge C_{1,\delta_0}\int\abso{\del u_0}^2$ and hence $E(u_0)>0$. Because of Proposition 5.3,  we only need to treat the case where there exists $\braces{t_n}_{n=1}^\infty$, $t_n\ge 0$, $t_n\uparrow T_+(u_0)$, so that 
$$\lambda(t_n)\to 0.$$  
(If $t_n\rightarrow t_0 \in [0,T_+(u_0))$, we obtain for all $R>0$, $\int_{|x| \ge R}|v(t_0)|^{2*}=0$ but  $\int\abso{\del v(t_0)}^2>0$).  After possibly redefining  $\braces{t_n}_{n=1}^\infty$ we can assume that 
$$
\lambda(t_n) \le 2 inf_{t\in [0,t_n]}\lambda(t)
$$  
and  from our hypothesis
$$w_n(x) = \frac{1}{\lambda(t_n)^{(N-2)/2}}u\paren{\frac{x}{\lambda(t_n)},t_n} \rightarrow w_0 \mbox{ \  in  \ }\dot H^1.$$  
By Theorem 3.9 we have $E(W)> E(w_0)=E(u_0)>0$, 
$\int\abso{\del u(t)}^2 \le (1-\overline\delta)\int\abso{\del W}^2$ so that $\int\abso{\del w_0}^2 < \int\abso{\del W}^2$.
 Thus $w_0\not\equiv 0$.  Let us now consider solutions of (CP), $w_n(x,\tau)$, $w_0(x,\tau)$ with data $w_n(-,0)$, $w_0(-,0)$ at $\tau = 0$, defined in maximal intervals $\tau\in (-T_-(w_n),0]$, $\tau\in (-T_-(w_0),0]$ respectively.  \\
Since $w_n\to w_0$ in $\dot H^1$, $\varliminf_{n\to \infty} T_-(w_n) \ge T_-(w_0)$ and 
$$
\mbox{ for each  \  } \tau\in (-T_-(w_0),0], \ \ w_n(x,\tau)\to w_0(x,\tau) \mbox{ \  in    } \dot H^1.  \mbox{  (See Remark 2.17)  } 
$$
Note that by uniqueness in (CP) (see Definition 2.10), for $0\le t_n + \tau/\lambda(t_n)^2$, $w_n(x,\tau) = \frac{1}{\lambda(t_n)^{(N-2)/2}} u\paren{\frac{x}{\lambda(t_n)},t_n+
\frac{\tau}{\lambda(t_n)^2}}$.  
Remark that $\varliminf_{n\to \infty} \tau_n = t_n\lambda(t_n)^2 \ge T_-(w_0)$ and thus for all $\tau\in (-T_-(w_0),0]$ for $n$ large, $0\le t_n + \tau/\lambda(t_n)^2 \le t_n$. Indeed,if $\tau_n \rightarrow \tau_0 < T_-(w_0)$, then  $w_n(x,-\tau_n) = \frac{1}{\lambda(t_n)^{(N-2)/2}} u_0(\frac{x}{\lambda(t_n)}) \rightarrow w_0(x,-\tau_0)$  in $\dot H^1$ with $\lambda(t_n)\to 0$ which is a contradiction from $u_0\not\equiv 0$, $w_0\not\equiv 0$.

Fix now $\tau\in (-T_-(w_0),0]$, for $n$ sufficiently large  $v(x,t_n+\tau/\lambda(t_n)^2)$, $\lambda(t_n+\tau/\lambda(t_n)^2)$ are defined and we have 

\begin{equation}\label{5.13}
\begin{aligned}
v(x,t_n+\tau/\lambda(t_n)^2) &= \frac{1}{\lambda(t_n+\tau/\lambda(t_n)^2)^{(N-2)/2}}u\paren{\frac{x}{\lambda(t_n+\tau/\lambda(t_n)^2)},t_n+\tau/\lambda(t_n)^2}\\
&= \frac{1}{\widetilde \lambda_n(\tau)^{\frac{N-2}{2}}}w_n\paren{\frac{x}{\widetilde \lambda_n(\tau)},\tau},
\end{aligned}
\end{equation}
with
\begin{equation}\label{5.14}
\widetilde \lambda_n(\tau) = \frac{\lambda(t_n+\tau/\lambda(t_n)^2)}{\lambda(t_n)} \ge \frac 1 2
\end{equation}
(because of the fact $\lambda(t_n) \le 2 inf_{t\in [0,t_n]}\lambda(t)$.)
One can assume after passing to a subsequence that $\widetilde \lambda_n(t_n+\tau/\lambda(t_n)^2) \rightarrow \widetilde \lambda_0(\tau)$ with $\frac 12 \le  \widetilde \lambda_0(\tau) \le + \infty$ and $v(x,t_n+\tau/\lambda(t_n)^2)\to v_0(x,\tau)$ in $\dot H^1$, as $n\to \infty$.
Remark that $\widetilde\lambda_0(\tau) < + \infty$. If not, we will have 
$\frac{1}{\widetilde \lambda_n(\tau)^{(N-2)/2}}w_0\paren{\frac{x}{\widetilde \lambda_n(\tau)},\tau} \rightarrow v_0(x,\tau)$
which implies $w_0(x,\tau) = 0$ which contradicts $E(w_0)=E(u_0)>0$. Thus $\widetilde \lambda_0(\tau) < + \infty$ and $v_0(x,\tau) = \frac{1}{\widetilde \lambda_0(\tau)^{(N-2)/2}}w_0\paren{\frac{x}{\widetilde \lambda_0(\tau)},\tau}$ where  $v_0(\tau) \in \overline K$. 
We thus obtain a contradiction from Proposition 5.3. Note that the same proof applies in the nonradial situation with the extra parameter $x(t_n)$.

\end{proof}

\begin{corollary}\label{5.15}
Assume that $E(u_0)<E(W)$, $\int\abso{\del u_0}^2<\int\abso{\del W}^2$ and $u_0$ is radial.  Then the solution $u$ of the Cauchy problem (CP) with data $u_0$ at $t=0$ has time interval of existence $I = (-\infty,+\infty)$, $\norm{S((-\infty,+\infty))}{u}<+\infty$ and there exists $u_{0,+}$, $u_{0,-}$ in $\dot H^1$ such that
\[
\lim_{t\to +\infty}\norm{\dot H^1}{u(t) - e^{it\Delta}u_{0,+}} = 0,\quad \lim_{t\to -\infty}\norm{\dot H^1}{u(t) - e^{it\Delta}u_{0,-}} = 0.
\]
Moreover, if we define $\delta_0$ so that $E(u_0) \le (1-\delta_0)E(W)$, there exists a function $M(\delta_0)$ so that $$\norm{S((-\infty,+\infty))}{u} \le M(\delta_0).$$
\end{corollary}
\begin{proof}
From the integral equation in Theorem 2.5, it is clear that $u(t)$ is radial for each $t\in I$.  Using Remark (4.23) and Theorem 5.1 we obtain (SC) or $I = (-\infty,+\infty)$, $\norm{S((-\infty,+\infty))}{u}<+\infty$.  Now Remark 2.15 finishes the proof of the first statement.\\
For the last statement, let
$$
D_{\delta_0} = \braces{ u_0 \in \dot H^1 \ \ \mbox{radial}, \ \ \int\abso{\del u_0}^2<\int\abso{\del W}^2 \ \mbox{ and } \ \ E(u_0) \le (1-\delta_0)E(W) }
$$
$$
M(\delta_0) = sup_{u \in D_{\delta_0}} \norm{S((-\infty,+\infty))}{u}.
$$
We need to show $M(\delta_0) < +\infty$. If not there is a sequence $u_{0,n}$ in $D_{\delta_0}$ and the corresponding solutions $u_n$ such that $\norm{S((-\infty,+\infty))}{u_n} \rightarrow +\infty$ as $n \rightarrow +\infty$. Note that we can assume that $\norm{S((-\infty,+\infty))}{e^{it\Delta}u_{0,n}} \ge \delta$, with  $\delta$ as in  Theorem 2.5. Arguing as in the proof of Proposition 4.1, we would conclude that first $J = 1$ in the decomposition given in Lemma 4.3 and then since $\norm{S((-\infty,+\infty))}{U_1}<+\infty$ we reach a contradiction (See also \cite{K}, Corollary 1.14).
\end{proof}

\begin{remark}\label{5.16}
Note that Corollary 5.14 is sharp.  In fact, $W(x)$ is radial and clearly $\norm{S((-\infty,+\infty))}{W} = +\infty$.  Moreover, Remark 3.14 shows that if $u_0\in H^1$ radial, $E(u_0)<E(W)$, but $\int\abso{\del u_0}^2>\int\abso{\del W}^2$, we have that $I$, the maximal interval of existence, is finite.
\end{remark}

Let us remark that we have, in fact, proved a slightly stronger result:

\begin{corollary}\label{5.16}
Let $u_0 \in \dot H^1$  be radially symetric and assume that for all \\
$t \in (-T_-(u_0),T_+(u_0))$  we have $\int\abso{\del u(t)}^2 \le \int\abso{\del W}^2 - \delta_0$, for $\delta_0>0$.  Then the solution $u$ of the Cauchy problem (CP) with data $u_0$ at $t=0$ has time interval of existence $I = (-\infty,+\infty)$, $\norm{S((-\infty,+\infty))}{u}<+\infty$.
\end{corollary}
\begin{proof}
Note first that if $E(u_0)<E(W)$, Corollary 5.14 yields the result, so that we can assume that $E(u_0) \ge E(W)$. Observe that the end of the proof of Lemma 3.4 gives us that, for $t \in (-T_-(u_0),T_+(u_0)),$
\begin{gather}\label{5.17}
\int \abso{\del u(t)}^2 - \int \abso{u(t)}^{2^*} \ge  C_{\delta_0} \int \abso{\del u(t)}^2.
\end{gather}
and that energy conservation, the assumption that $E(u_0) \ge E(W)$ yields $$inf_{t \in (-T_-(u_0),T_+(u_0))} \int \abso{\del u(t)}^2  \ge C.$$

From the Remark 2.7, if $\delta_o$ is close to $\int\abso{\del W}^2$, our conclusion holds. We can then find $0 \le \delta_c < \int\abso{\del W}^2$, so that if for $t \in (-T_-(u_0),T_+(u_0)), \int\abso{\del u(t)}^2 \le \int\abso{\del W}^2 - \delta_0, \delta_0 > \delta_c,$ our desired conclusion holds and $\delta_c$ is optimal with this property. We assume $\delta_c>0$ and reach a contradiction by establishing the analogs of Propositions 4.1,  4.2 and using (5.17). In order to do so, we just need to modify the statement of Lemma 4.9, by replacing (4.10), say, by for $t \in I$,
$\int\abso{\del U_1(t)}^2 \le \int\abso{\del W}^2 - \delta_c$, where $U_1(t)$ is the non-linear profile in Lemma 4.9 and $I$ its maximal interval of existence as in Remark 2.13. (4.11) is replaced analogously. The proofs are then identical to those given in section 4 and that of Corollary 5.14.
\end{proof}

As a consequence of Corollaries 5.14 and 5.16, we obtain the following concentration phenomenon for all radial type II finite blow-up solutions. Here by type II finite blow-up solution, we mean a solution $u$ whose maximal interval of existence $I$ is finite and for which there is a $C$, such that for all $t\in I$, $\int \abso{\del u(t)}^2 < C$. On the other hand, type $I$ finite blow-up solution is such that what the time of existence is finite but the $ \dot H^1$ norm blows up.

\begin{corollary}\label{5.18}
Let $u_0 \in \dot H^1$ (no size restrictions) be radially symetric and assume that $T_+(u_0) < + \infty$, and that $\forall t \in[0,T_+(u_0)), \int\abso{\del u(t)}^2 \le C_0$. Then, for all $R>0$, we have
$$\varliminf_{t\to T_+(u_0)} \int_{|x| \le R}  \abso{\del u(t)}^2 \ge \frac 2 N \int  \abso{\del W}^2$$
$$\varlimsup_{t\to T_+(u_0)} \int_{|x| \le R}  \abso{\del u(t)}^2 \ge \int  \abso{\del W}^2.$$
\end{corollary}
\begin{proof}
Consider $t_n \to T_+(u_0)$ and apply Lemma 4.3 to the sequence $u(t_n)$. Arguing in an analogous manner to the proof of Theorem 2.14, we must have $\lambda_{j,n} \to +\infty$ for some $j$ and the corresponding non-linear profile $U_j$ has $\norm{S((0,T_+(U_j))}{U_j}=+\infty$. 
If the first inequality does not hold, we can find a sequence $\{t_n\}$ as before and $R_0>0$, $\eta_0>0$ so that 
$$ \int_{|x| \le R_0}  \abso{\del u(t_n)}^2 \le \frac 2 N \int  \abso{\del W}^2 - \eta_0.$$
In addition,we must have (since $\lambda_{j,n} \to +\infty$) that 
$$
 \int \abso{\del U_j(-t_{j,n}/\lambda_{j,n}^2)}^2 \le \frac 2 N \int  \abso{\del W}^2 - \eta_0 < \frac 2 N \int  \abso{\del W}^2 < \int  \abso{\del W}^2,
$$
Thus $E(U_j) < E(W) =  \frac 1 N \int  \abso{\del W}^2 $ and Corollary 5.14 gives a contradiction.\\
If the second inequality does not hold, we can find $R_0>0$, $\eta_0>0$ so that for all $t \in I$, $ \int_{|x| \le R_0}  \abso{\del u(t)}^2 \le \int  \abso{\del W}^2 - \eta_0$.
By the argument at the begining of the proof of case 1 of Proposition 5.13, we must have $-t_{j,n}/\lambda_{j,n}^2 < C$. Thus, we obtain, for $t>M$, that $\int  \abso{\del U_j(t)}^2 \le \int  \abso{\del W}^2 - \eta_0$, so that Corollary 5.16 concludes the proof.
\end{proof}

\begin{remark}\label{5.19} Note that we have not yet shown that $u_0$ as in Lemma 5.18 exist, but we expect that this is the case. We also expect to have data $u_0$ for which type I blow-up 
 occurs.
\end{remark}

\begin{remark}\label{5.20} In the case $N \ge 4$, consider now $u_0 \in H^1$ radial as in Corollary 5.18 (but not type II), then using the $L^2$ conservation and energy laws, estimates as in \cite{MT} yield for any sequence $t_n$ such that $\int \abso{\del u(t_n)}^2 \to +\infty$ that for all $R>0$, we have $\int_{|x| \le R}  \abso{\del u(t_n)}^2 \to +\infty$ which leads to the same conclusions as in Corollary 5.18. Note that when $N=3$, one expects that the conclusion in this remark is false in light of examples analogous to the ones constructed by Raphael in \cite{R} which give a radial solution blowing up exactly on a sphere.
\end{remark}

\nocite{*}


\begin{thebibliography}{10}


\bibitem{A}   Aubin, Thierry,    \'{E}quations diff\'erentielles non lin\'eaires et probl\`eme              de {Y}amabe concernant la courbure scalaire,  J. Math. Pures Appl. (9),    55,      1976, 3, 269--296


\bibitem{BG} Bahouri, Hajer; G{\'e}rard, Patrick High frequency approximation of solutions to critical nonlinear wave equations.  Amer. J. Math. 121 (1999), no. 1, 131--175. 

\bibitem{BC}  Berestycki, Henri; Cazenave, Thierry,
Instabilit{\'e} des {\'e}tats stationnaires dans les {\'e}quations de  Schr\"odinger et de Klein-Gordon non lin{\'e}aires
C. R. Acad. Sci. Paris SŽr. I Math. 293 (1981), no. 9, 489--492.


 \bibitem{BL}  Bergh, Joran; Lofstrom, Jorgen Interpolation spaces. An introduction. Grundlehren der Mathematischen Wissenschaften, No. 223. Springer-Verlag, Berlin-New York, 1976.


\bibitem{B1}    Bourgain, J.,   Global well-posedness of defocusing critical nonlinear              {S}chr\"odinger equation in the radial case,   J. Amer. Math. Soc.,    12,     1999,  1,   145--171

\bibitem{B2}   Bourgain, J.,   New global well-posedness results for nonlinear Schr\"odinger equations,   AMS Publications,     1999


\bibitem{C}   Cazenave, Thierry,   Semilinear {S}chr\"odinger equations,    Courant Lecture Notes in Mathematics,    10, New York University Courant Institute of Mathematical              Sciences,   New York,    2003

\bibitem{CW}   Cazenave, Thierry and Weissler, Fred B.,    The {C}auchy problem for the critical nonlinear              {S}chr\"odinger equation in {$H\sp s$},   Nonlinear Anal.,  Theory, Methods \&  Applications. An              International Multidisciplinary Journal. Series A: Theory and              Methods,    14,      1990,   10,    807--836


\bibitem{CKSTT}   Colliander, J. and Keel, M and Staffilani, G. and Takaoke, H. and Tao, T.,    Global well-posedness and scattering for the energy-critical nonlinear Schr\"odinger equation in {$\mathbb{R}^3$},   to appear, Annals of Math


\bibitem{Ge} G{\'e}rard, Patrick Description du dŽfaut de compacit{\'e} de l'injection de Sobolev ESAIM Control Optim. Calc. Var. 3 (1998), 213--233


\bibitem{GMO}   Gerard, Patrick and Meyer, Yves and Oru, Fr{\'e}d{\'e}ric,   In\'egalit\'es de {S}obolev pr\'ecis\'ees, S\'eminaire sur les \'Equations aux D\'eriv\'ees Partielles,              1996--1997, Exp.\ No.\ IV, 11, \'Ecole Polytech.,   1997

\bibitem{G}   Glassey, R. T.,    On the blowing up of solutions to the {C}auchy problem for              nonlinear {S}chr\"odinger equations,   J. Math. Phys.,   18,      1977,   9,    1794--1797


\bibitem{Gr}   Grillakis, Manoussos G.,    On nonlinear {S}chr\"odinger equations,   Comm. Partial Differential Equations,    25,     2000,    9-10,     P1827--1844


\bibitem{KT}   Keel, Markus and Tao, Terence,    Endpoint {S}trichartz estimates,   Amer. J. Math.,    120,      1998,    5,    955--980


\bibitem{K}   Keraani, Sahbi,   On the defect of compactness for the {S}trichartz estimates of              the {S}chr\"odinger equations,   J. Differential Equations,    175,      2001,    2,     353--392

\bibitem{K1}   Keraani, Sahbi, On the blow up phenomenon of the critical {S}chr\"odinger equation, J. F. A., V235 (2006), 171--192


\bibitem{M}   Merle, F.,  Determination of blow-up solutions with minimal mass for              nonlinear {S}chr\"odinger equations with critical power,   Duke Math. J.,    V69,   1993,    2,   427--454

\bibitem{M1} Merle, F., On uniqueness and continuation properties after blow-up time of self-similar solutions of nonlinear Schrödinger equation with critical exponent and critical mass.  Comm. Pure Appl. Math.  45  (1992),  no. 2, 203--254

\bibitem{M2} Merle, Frank,Existence of blow-up solutions in the energy space for the critical generalized KdV equation.  J. Amer. Math. Soc.  14  (2001),  no. 3, 555--578

\bibitem{MT} Merle, Frank and Tsutsumi, Yoshio, $L\sp 2$ concentration of blow-up solutions for the nonlinear Schrödinger equation with critical power nonlinearity.  J. Differential Equations  84  (1990),  no. 2, 205--214. 


\bibitem{MV} Merle, F.; Vega, L. Compactness at blow-up time for $L\sp 2$ solutions of the critical  nonlinear Schr\"odinger equation in 2D.  Internat. Math. Res. Notices 1998, no. 8, 399--425.


\bibitem{OT}  Ogawa, Takayoshi; Tsutsumi, Yoshio Blow-up of $H\sp 1$ solution for the nonlinear Schr\"odinger equation.  J. Differential Equations  92  (1991),  no. 2, 317--330

\bibitem{R}   Raphael, Pierre, Existence and stability of a solution blowing-up on a sphere for a $L^2$ supercritical nonlinear Schrodinger equation, Duke Math. J. (to appear)

\bibitem{RV}   Ryckman, E. and Visan,  Monica,    Global well-posedness and scattering for the defocusing energy-critical nonlinear Schrodinger equation in {$\mathbb{R}^{1+4}$},   preprint, http://arkiv.org/abs/math.AP/0501462

\bibitem{Ta} Talenti,  Giorgio, Best constant in Sobolev inequality,  Ann. Mat. Pura Appl., (4)  110  1976, 353--372

\bibitem{T}   Tao, Terence,   Global well-posedness and scattering for the              higher-dimensional energy-critical nonlinear {S}chr\"odinger              equation for radial data,   New York J. Math.,    11,      2005,    57--80 (electronic)


\bibitem{TV}   Tao, Terence and Visan, Monica,   Stability of energy-critical nonlinear {S}chr\"odinger              equations in high dimensions,   Electron. J. Differential Equations,    2005,      118, 28 pp. (electronic)



\bibitem{V} Visan, M.,   The defocusing energy-critical nonlinear Schr\"odinger equation in higher dimensions, preprint, http://arkiv.org/abs/math.AP/0508298




\bibitem{W} Weinstein, Michael
Nonlinear Schr\"odinger equations and sharp interpolation  estimates. 
Comm. Math. Phys. 87 (1982/83), no. 4, 567--576.


\bibitem{Z}  Zhang, Jian Sharp conditions of global existence for nonlinear Schr\"odinger and  Klein-Gordon equations. Nonlinear Anal. 48 (2002), no. 2, Ser. A: Theory Methods, 191--207.

\end{thebibliography}
\end{document}